\input amstex
\magnification=1200
\loadmsam
\loadmsbm
\loadeufm
\loadeusm
\UseAMSsymbols

\hsize=6.9truein
\hoffset=-0.11truein
\vsize=8.9truein
\voffset=-0.2truein

\def\leftitem#1{\item{\hbox to\parindent{\enspace#1\hfill}}}

\def\boxit#1#2{\hbox{\vrule
	\vtop{%
	\vbox{\hrule\kern#1%
	\hbox{\kern#1#2\kern#1}}%
	\kern#1\hrule}%
	\vrule}}

\def\leaderfill{\leaders\hbox to 1em{\hss.\hss}\hfill}

\parskip=\medskipamount
\document

\input epsf



\centerline{\bf Grothendieck's Reconstruction Principle and}
\centerline{\bf 2-dimensional Topology and Geometry}

\medskip
\centerline{\bf Feng Luo}

\S0. {\bf Introduction}

The goal of this paper is an attempt to relate some ideas of Grothendieck
in his Esquisse d'un programme [Gr1] and some of the recent results  
on 2-dimensional topology and geometry. Especially, we shall  discuss the Teichm\"uller
theory, the mapping class groups, $SL(2, \bold C)$ representation variety
of surface groups, and Thurston's theory of measured laminations.

A prominent idea in surface theory  is that to study a surface, one should consider
all subsurfaces inside it. 
Indeed, there is a 
\it hierarchy \rm
of compact oriented surfaces of negative Euler number under (essential)
inclusion. Each surface in the hierarchy is indexed by its \it level \rm
which is the number of disjoint simple loops needed to decompose it into
3-holed spheres (i.e., the complex dimension of the Teichm\"uller space 
of complete hyperbolic metrics of finite area). The first three levels
in the hierarchy are listed as follows. The level-0 surface is the  3-holed
sphere, the level-1  surfaces are 1-holed torus and 4-holed sphere, and
the level-2 surfaces are 2-holed torus and 5-holed sphere.

\midspace{0.1cm}
\centerline{\epsfbox{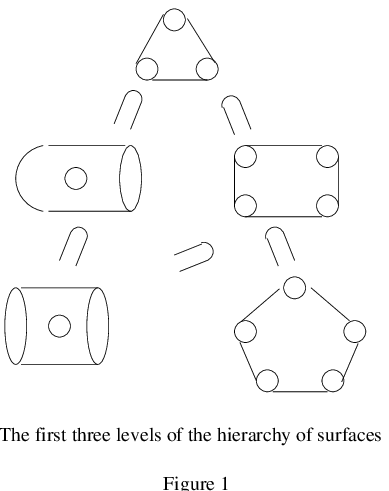}}

One of the key ideas in [Gr1] which we would like to discuss at
 length in this paper is Grothendieck's reconstruction principle
for the ``Teichm\"uller tower". We quote the relevant paragraph
in [Gr1] below. On page 11, line 2-6, Grothendieck wrote (with
English translation by L. Schneps)
`` The a priori interest of a complete knowledge of the first
two levels of the tower is to be found in the principle that
$\underline{\text{
the
entire tower can be reconstructed from these two first levels}}$,
in the sense that via the fundamental operation of `gluing',
level-1 gives the complete system of generators, and level-2
a complete system of relations."
One may interpret this  principle broadly as follows.
To study  a structure 
(for instance, hyperbolic structure, complex projective
structure, measured lamination, or linear representation of the 
surface group) and its moduli space on a surface, 
one should consider the restrictions of the
structure to the level-1 subsurfaces and reconstruct the structure 
from its restrictions.
The level-2 surfaces should serve as ``relators" in the reconstruction
process. For example, one may ask if the reconstruction principle holds
for the characters of representations of the surfaces groups into 
 general linear group $GL(n, \bold C)$.  Namely, suppose $f$ is a complex
valued function defined on the fundamental
group of the surface so that the restriction of $f$ to
the fundamental group of
each essential level-1  subsurface is a $GL(n, \bold C)$-character. Is $f$
the character of some $GL(n, \bold C)$  representation of the surface group?
In [Lu3], we show that the answer is affirmative for $SL(2, \bold C)$
representations of the surface groups.
It is interesting to note that this principle 
of reconstruction was taken as one 
of the basic axioms by physicists in conformal field theory ([MS]).

The main theorems in [Lu1], [Lu2] state that the Teichm\"uller
space and Thurston's measured lamination space
for surfaces obey the reconstruction principle. Also using the
work of Gervais [Ge], we see that
the mapping class group of a surface fits the principle as well [Lu4].
These will be the main topics of this paper.
We shall  also discuss
some open questions arising from reading of [Gr1].


We remark that as far as we know, there is no precise definition of
the Teichm\"uller tower in  [Gr1]. See also the books  [Gr2],
and [Gr3].
What follows is my interpretation of
Grothendieck's reconstruction principle and there should be other ways
of interpreting  it (for instance in algebraic geometry).

One way to illustrate Grothendieck's reconstruction principle is to
consider convex planar $n$-sided polygons. According to the
principle, to construct a convex $n$-sided ($n \geq 5$) polygon, one should
consider all convex quadrilaterals inside the polygon (the vertex of the
quadrilateral is a vertex of the polygon).  The convex polygon
is a union  of these quadrilaterals by gluing along their overlaps. 
Now these quadrilaterals overlap in two different ways. An  \it essential
overlap \rm of two quadrilaterals contains an edge or diagonal. Otherwise,
they overlap inessentially (see figure 2). The reconstruction principle
states that it suffices to glue quadrilaterals along the essential overlaps.
The gluing along the inessential  overlaps is a \it consequence \rm of
the gluing along essential overlaps (see \S1 for more details). As a consequence,
to study the geometry of the moduli space of convex  polygons, it suffices
to understand that of quadrilaterals.

\midspace{0.1cm}
\centerline{\epsfbox{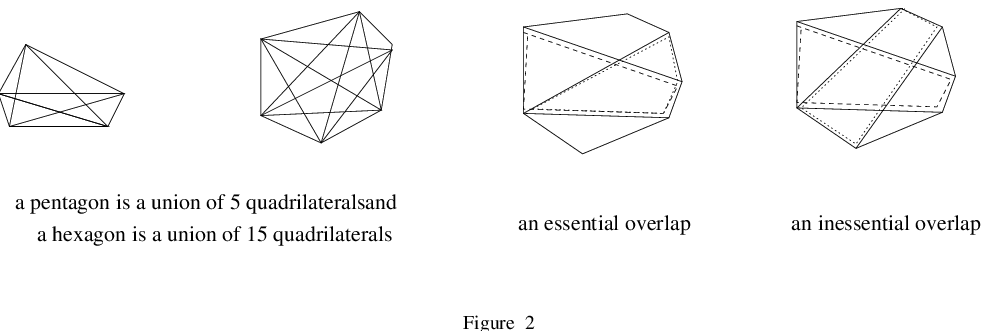}}

The situation for
surfaces is analogous to that of polygons where the 3-holed sphere 
corresponds to the  triangle and 1-holed torus and 4-holed sphere correspond
to quadrilaterals. (One should think of the polygons forming an hierarchy
under inclusion. And the level of a polygon is the number of disjoint
diagonals needed to decompose a polygon into triangles.) 
Thus according to the reconstruction principle, to construct
a hyperbolic metric on a surface of negative Euler number, one should
consider all (isotopy classes of) subsurfaces which are homeomorphic to
the 1-holed torus or the 4-holed sphere. These subsurfaces overlap in two
different patterns: the overlap  is \it essential \rm if there is a
homotopically non-trivial loop in the overlap, otherwise it is inessential
(see figure 3). 
The reconstruction principle says that we can glue along essential
overlaps  to recover the original hyperbolic structure. 
To be more precise,
assign to each level-1 subsurface  a hyperbolic structure so that
when two level-1 subsurfaces overlap essentially, they overlap geometrically
(i.e., the geodesic lengths of  all overlapping simple loops are the same
in both level-1 surfaces).  Then the
reconstruction principle states that there exists
a hyperbolic metric on the surface whose restrictions to level-1
subsurfaces are  (isotopic to) the assigned hyperbolic structures.
This is the main result established in [Lu1]. 

\midspace{0.1cm}
\centerline{\epsfbox{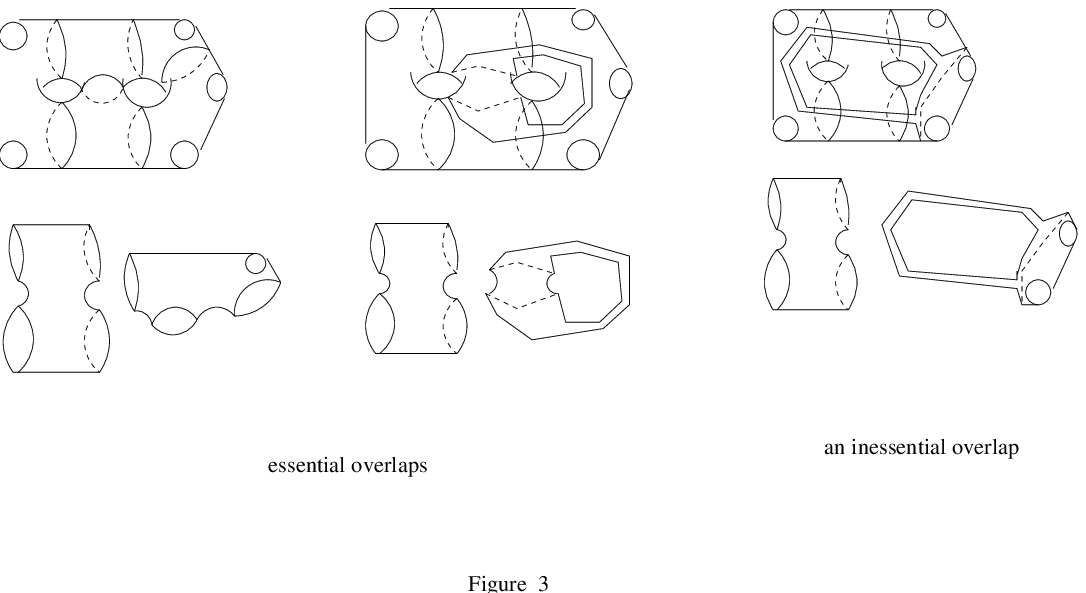}}

The organization of the  paper is as follows. In \S1, we study the
moduli space of convex polygons in details and use it to
illustrate Grothendieck's principle. We also discuss the  ideal triangulations of surfaces.
In \S2, we recall the Teichm\"uller spaces, the mapping class groups,
and Thurston's projective measured lamination spaces. 
It is well known that these
three themes represent the geometric, algebraic, and topological aspects
of surface theory.  In the case of the torus, these three themes correspond
to the upper-half space $\bold H$, $SL(2, \bold Z)$, and $\hat  \bold R  = \bold R\cup
\{\infty\}$ which  appeares as the boundary of $\bold H$. The natural action of the mapping
class group corresponds to the action of $SL(2,\bold Z)$ on $\bold H$ and $\bold
R$ by M\"obius transformations.  We shall also recall the related topics for
level-1  surfaces.
In \S3, we state the reconstruction theorems for the Teichm\"uller spaces,
the measured lamination spaces and the mapping class groups.  
In \S4,  we discuss the key ingredient
in the proofs of the reconstruction theorem, namely simple loops on
surfaces. 
We also recall the notion of  $SL(2, \bold Z)$
modular structures on a set.  The role of the modular structure, or
the same, $(\bold Q P^1, SL(2, \bold Z))$ is prominent in the 
reconstruction program as predicted by Grothendieck (see page 
248-249 in [Gr1] or \S 3.8).
Topologists have been knowing the role of  modular configuration for
simple loops on level-1 surfaces since the fundamental work of Max
Dehn [De] in 1938. Dehn actually used the structure to give an
elegant derivation of the mapping class group of the  4-holed sphere
(see \S 4.1).  The special feature of the modular configuration
is the huge symmetry built in the configuration. This is, in our view,
one reason why the set of homotopy class of simple loops on the surface
is more useful than the fundamental group in establishing 
the reconstruction principle for many structures (see \S 3.7, figure 8,
and \S 6.1).  
In \S5 we give a fairly general reason
which indicates the special role played by level-2 surfaces in the
reconstruction principle.
In the last section, we  discuss the characters of
$SL(2, \bold C)$ representations.

This is not a general survey of recent advances on surface theory.  And
we have omitted the important contributions of Masur and Minsky and others
on 2-dimensional topology and geometry because of our lack of expertise in
the areas. 

{\it Acknowledgement.} I would like to thank X.S. Lin for inviting
me to write the paper.  Discussions with F. Bonahon, 
L. Keen, X.-S. Lin and C. Series have been very helpful for
me in developing ideas in the paper. I thank the referee and
P. Landweber for careful reading of the manuscript and for
suggestions on improving the exposition of the paper. 
The work is supported in part by the NSF.


\S1. {\bf A Simple Example of Convex Polygons}

1.1.
We shall illustrate the reconstruction principle and its
application by considering the configuration space of convex $n$-sided polygons.
Let us begin with the following problem.

\noindent
{\bf Problem 1.} Describe the space $T(n)$ of all convex $n$-sided polygons
up to isometries. Here polygons have marked vertices and isometries preserve
markings.

To be more precise, let us distinguish the topological (or combinatorical)
and geometric aspects of the problem. By an $n$-sided polygon we mean a
topological disk with $n$ marked points (the vertices) on its boundary.
A \it convex structure \rm on a polygon is a metric on the polygon which
is isometric to a convex $n$-sided planar polygon so that the marked points
correspond to the vertices.  An \it edge \rm in a convex polygon is a
line segment joining two vertices, and a \it diagonal \rm is an
edge which does not joint adjacent vertices. Given an $n$-sided polygon $P$, the space
of all convex structures on $P$ modulo isometries preserving the
vertices is denoted by $T(P)$ which is essentially $T(n)$. 

Having introduced these notations, we may rephrase the problem as follows.

\noindent
{\bf Problem 2.} Assign to each edge in a convex $n$-sided polygon a 
positive number.  When is the assignment 
corresponds to the edge lengths of an $n$-sided convex polygon?

The solution for triangles $n=3$ is well known. Namely, the
assignment must satisfy the triangular inequalities that sum of
two is larger than the third. For general $n$,
the assignment must satisfy triangular inequalities over three
edges forming a triangle and equations over six edges forming a
quadrilateral. Grothendieck's reconstruction principle asserts that
these are the set of all constrains, i.e., the quadrilaterals 
(= level-1 polygon) are the ``generators" in building convex polygons.

1.2.  It is instructive to compare the hierarchies of polygons
and surfaces.
One first observes that the isometry class of a convex polygon
is determined by the lengths of all edges. 
The corresponding fact in hyperbolic geometry is a result of
Fricke-Klein [FK] that the isometry class of a hyperbolic metric on
a surface is determined  by the lengths of simple geodesic loops.
The solution for $T(3)$ is given by  $T(3) =\{(a_1, a_2, a_3)
\in \bold R^3 | a_i + a_j > a_k\}$ reflecting the fact that a triangle
is determined up to isometry by its three edge lengths subject to the
triangular inequalities. The corresponding fact in hyperbolic geometry
is the well known theorem of Fricke-Klein [FK] that a hyperbolic
metric on a 3-holed sphere is determined up to isometry by the three
lengths of the boundary geodesics and these  lengths subject no constraints.
For $n \geq 4$, an old way of solving the problem for $T(n)$ is to
triangulate the $n$-sided polygon by $(n-3)$ edges,  i.e., one uses the
triangle as the basic building block. This corresponds to the Fenchel-Nielsen's
decomposition of surfaces into 3-holed spheres (using 3-holed sphere as
the basic building block). In this way, one parametrizes the convex polygon
by the lengths of the edges in the triangulation.
These lengths have to satisfy complicated  inequalities due to the convexity.
Unlike the Fenchel-Nielsen coordinate for Teichm\"uller space which can be used
to express the Weil-Petersson symplectic form  by Wolpert's formula
([Wo]), the length coordinate
for convex polygons  seems to be  
less useful in extracting geometric information
about $T(n)$ except that it can be used to show that $T(n)$ is a real analytic
manifold diffeomorphic to $\bold R^{2n-3}$.

1.3.
Now Grothendieck's reconstruction principle asserts that quadrilateral be the
basic building block. 
Thus given an $n$-sided polygon $P$, one considers all
quadrilateral whose edges are edges in $P$. If the polygon
has a convex structure, each quadrilateral in $P$ becomes a convex
quadrilateral. These convex quadrilaterals satisfy the obvious consistency
condition:

\noindent (*)  \quad
If two quadrilaterals overlap  essentially (i.e., there is  an edge 
in the overlap), then the convex quadrilaterals overlap 
geometrically, i.e., the corresponding lengths of edges 
in both convex quadrilaterals are the same. 

It turns out that the condition (*) is also sufficient to recover
the convex polygon for obvious reason.

\noindent
{\bf Reconstruction Principle for Polygons}. \it To construct a  convex
$n$-sided polygon with $n \geq 5$, it suffices to assign to each quadrilateral
in the polygon a convex structure so that the assignment satisfies the
consistency condition (*). \rm

1.4. In terms of the reconstruction principle, the solution to
problem 2 is simply that the assignment must be realized by
a convex quadrilateral  for each choice of six edges forming
a quadrilateral.

This principle  essentially reduces the study of $T(n)$ to that of 
$T(4)$. To understand $T(4)$, one uses the lengths of the  six
edges of a quadrilateral. First of all, the lengths satisfy the
triangular inequalities that sum of two is larger than the third
 over each of the four triangles in the
quadrilateral.
By a simple calculation, one shows that
these six lengths
satisfy the following constraint:

\midspace{0.1cm}
\centerline{\epsfbox{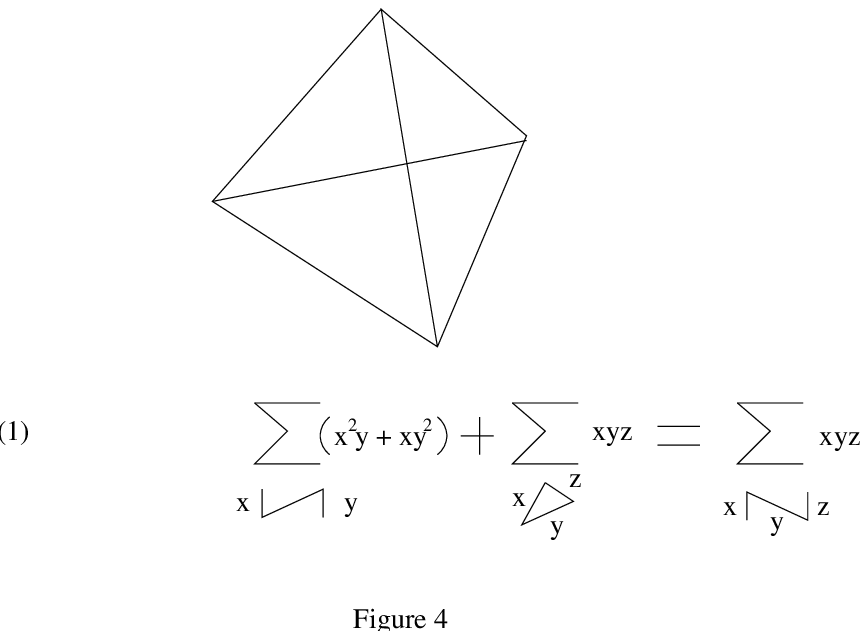}}
\midspace{0.1cm}

\noindent 
where the sums are over all specified subgraphs (of the complete graph
on 4 vertices) whose edges
are labeled by their lengths.
The convexity condition is equivalent to the \it Largest Root Condition
\rm  that the lengths of the
diagonals are the largest root. (Fix five lengths and think of (1) as
a quadratic equation in the length of the remaining  diagonal. It has two real roots and
the convexity condition says that the largest root is the length.)
As a consequence of (1), the Largest Root Condition, the triangular
inequalities,
and the reconstruction principle, one obtains 
a complete solution to problem 2. 

1.5.
Given a convex polygon $P$, the observable  invariants of $P$ seem to be the
area of the polygon, the lengths of edges and the angles of 
intersections of edges. These
define  the ``observable" area, length and angle functions on 
the configuration space
$T(P)$. To be more precise, fix an isotopy class of an edge  
$e$ in $P$ (resp. a pair of isotopy classes of intersecting edges), one defines a length function (resp. angle function) from $T(P)$ to $\bold R$ by sending a convex structure to the
length  (resp. angle) of  $e$ in the convex structure.
These naturally defined functions seem to play an important role
in the geometry of the configuration space $T(P)$. 
And indeed they are. Here is one way to see it using Thurston's
invariant of oriented triangles.
Suppose $\Delta ABC$ is an oriented triangle in the plane so that
the cyclic order $(A,B,C)$ is the right-hand orientation in the plane
(for simplicity, we shall assume triangles are right-hand oriented
in the plane unless stated otherwise). Then  the \it Thurston invariant \rm
$z_A$ of the triangle at edge $BC$ is defined to be the complex number
$\frac{C-A}{B-A}$ in the upper-half plane $\bold H$.
The edge invariants $z_B = \frac{A-B}{C-B}$ and $z_C = \frac{B-C}{A-C}$ 
for $AC$ and $AB$ are
given by $z_B = 1/(1-z_A)$ and $z_C = 1/(1-z_B) = (z_A -1)/z_A$ respectively. 
In particular $z_A z_B z_C = -1$. Evidently if two oriented triangles differ
by a similarity transformation ($ f(w)= aw + b, $ $a, b \in \bold C, a \neq 0$),
then their Thurston invariants are the same. For a convex quadrilateral
with a marked diagonal, one defines the Thurston invariant to  be the  pair
$(z, w) \in \bold H \times \bold H$ where each coordinate is the invariant
of a triangle at the marked diagonal. For instance, in terms of
Thurston invariants, parallelograms are exactly those convex quadrilaterals
with Thurston invariant $(z,z)$.
A simple calculation shows that
if the marked diagonal is changed, then the invariant becomes
$(u,v)$ where $ u = \frac{w-1}{w(1-z)}$ and
$v = \frac{z-1}{(1-w)z}$ as shown in figure 5. 

\midspace{0.1cm}
\centerline{\epsfbox{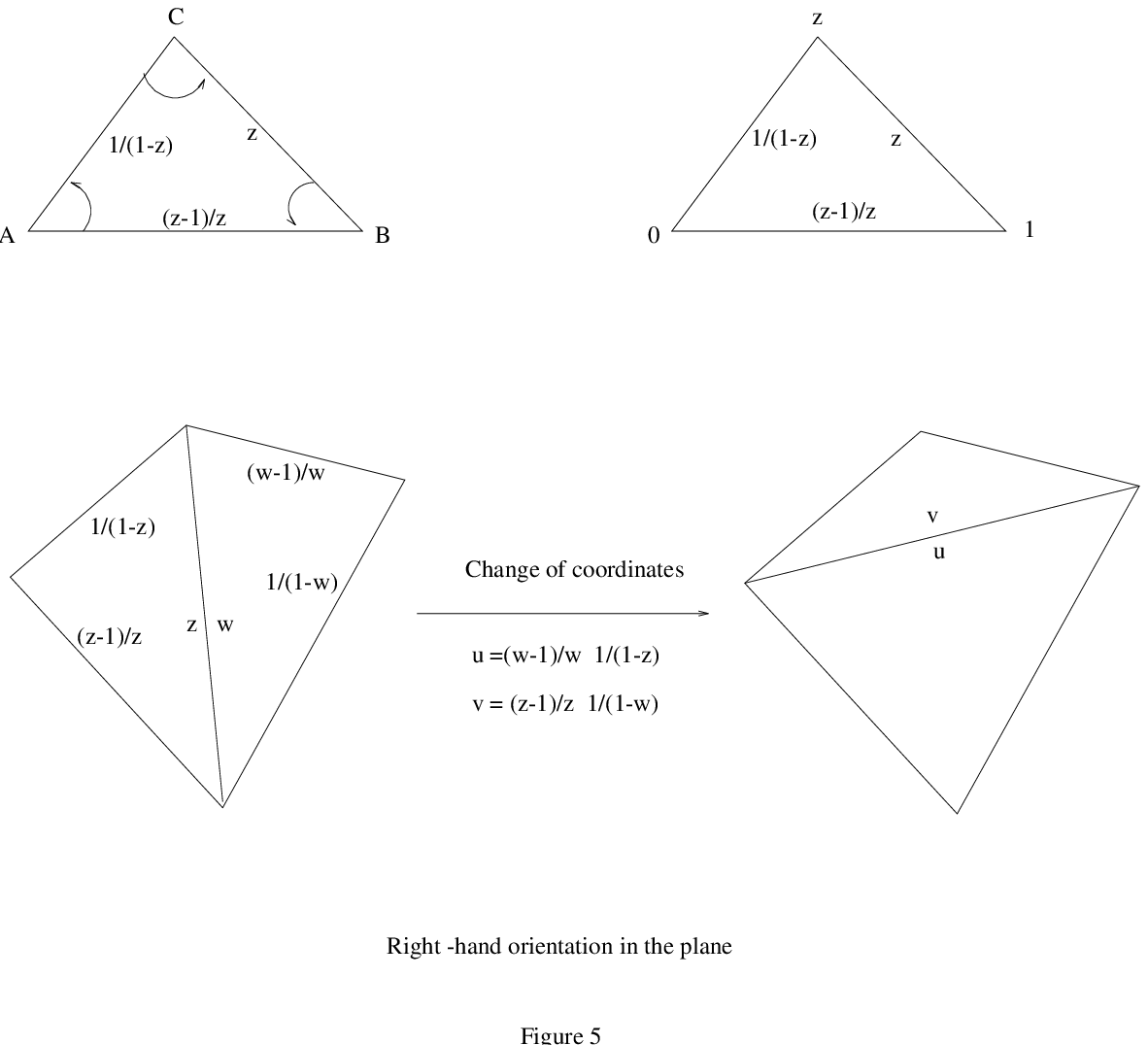}}

\noindent
Furthermore the convexity condition is
equivalent to both $(z,w)$ and $(u,v) \in \bold H \times \bold H$. As a consequence of
the transformation formula, one sees that the
space of similarity classes of convex quadrilaterals has a natural
complex structure.  Combining these with
the reconstruction principle, one obtains the fact that the projectivized
space $T(n)/\bold R_+$ of similarity classes of convex polygons has
a natural complex structure so that
angle functions and the logarithm of the ratio of the length 
functions are pluriharmonic. Furthermore, the space $T(n)/\bold R_+$
can be explicitly described. This result itself is not surprising
since another way of parametrizing $T(n)/\bold R_+$ is by taking
vertices as the  coordinates. But the fact that the complex structure
is built on that of $T(4)/\bold R_+$ seems to be interesting.
Evidently these ``observable" length and angle functions also
exist on the Teichm\"uller spaces. It  is natural to ask if these
functions are somehow related to the complex structure on the
Teichm\"uller space.

1.6. A close surface of genus $g $ is obtained by identifying
the opposite edges of a $4g$-sided polygon. Suppose that the
genus $g > 1$. Now take a 
$4g$-sided polygon bounded by line segments
and identify the opposite sides by affine
transformations. 
The quotient surface has a singular affine structure 
with one singularity coming from the vertices
(the monodromy around the singularity is of the form 
$z \to k z + (k-1)a$, $k \in \bold R_{>0}$). If immersed
$4g$-sided polygons are allowed to be  used in the construction, then all
complex structures on the surface are obtained in this way.
See [Lu7] for more details on relating the
Thurston's coordinate, the singular affine structure and
the complex structure on the Teichm\"uller space. 

1.7.
The reconstruction principle also  holds for hyperbolic or spherical
convex polygons. Thus the same picture holds in these cases as well.
The most interesting one  seems to be the ideal polygons in hyperbolic
plane where one assigns to each oriented ideal quadrilateral with
a marked diagonal the Bonahon-Thurston shearing coordinate ([Bo1], [Th1]).
Recall that the shearing coordinate is defined as follows. 
The \it mid-point \rm of
an edge  in an ideal triangle is the tangent point of the inscribed circle
at the edge. Given an oriented ideal quadrilateral with a marked diagonal,
one chooses an orientation on the diagonal. The Bonahon-Thurston
 coordinate for the marked quadrilateral is the exponential of the 
signed hyperbolic distance from the left mid-point to the right mid-point 
of the diagonal.
Note that the coordinate is independent of the choice of the orientation on the
diagonal.  If one changes the diagonal, the coordinate changes to its
inverse.
The change of coordinate formula for other four edges is 
given in  figure 6.

\midspace{0.1cm}
\centerline{\epsfbox{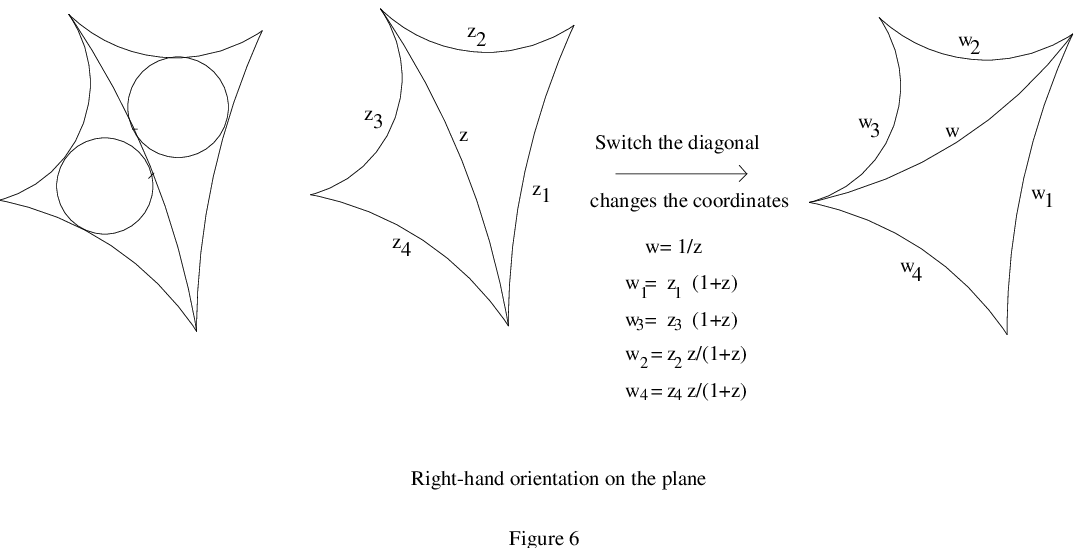}}

\noindent
Note that the transformation formulas  are real algebraic. Since 
each non-close surface has
an ideal triangulation, this gives an easy way to parametrize the
Teichm\"uller space of the surface using Bonahon-Thurston coordinates.
As a consequence, one proves easily that the Teichm\"uller space is
a real analytic manifold diffeomorphic to the Euclidean space.
(This seems to be one of the quickest way of showing that the Teichm\"uller
space of a non-closed surface is contractable. The other proof using the
Fenchel-Nielsen coordinate  seems to be always running into the technical
difficulties of  showing that  Fenchel-Nielsen coordinates
associated to  different 3-holed sphere decompositions differ
by a diffeomorphisms.)  These shearing coordinates are closely related
to Penner's coordinates for decorated Teichm\"uller spaces (see [Pe]).
See also [Mo] for related material on measured laminations.

\S2. {\bf The  Teichm\"uller Space, the Mapping Class Group, and 
the Space of Measured Laminations}
 
2.1. Given a compact orientable surface $\Sigma$ with or without boundary,
there are three themes naturally associated to the surface. Namely, the
Teichm\"uller space $T(\Sigma)$, the mapping class group  $\Gamma(\Sigma)$, and the
space  $\Cal S(\Sigma)$ of isotopy classes of unoriented
simple loops not homotopic to a point
(or its completion, Thurston's space of measured laminations).
These three themes represent the geometric, algebraic, and topological
aspects of the surface theory. Recall that the mapping class group $\Gamma(\Sigma)$
= Homeo$^+(\Sigma, \partial \Sigma)/$Iso is the  group of orientation
preserving self-homeomorphisms
modulo isotopies so that the boundary of the surface is fixed pointwise
by the homeomorphisms and the isotopies. For surface of negative Euler
number, the Teichm\"uller space $T(\Sigma)$ is the space of
all hyperbolic
metrics with geodesic boundary on the surface modulo isometries isotopic
to the identity. These three themes
interact each other in the sense that the mapping class group acts naturally
on both  $T(\Sigma)$ and $\Cal S(\Sigma)$ by pull back, and the space 
$\Cal S(\Sigma)$ appeared in Thurston's
compactification of the Teichm\"uller space $T(\Sigma)$.

2.2. One may illustrate these three themes and their interaction by the
classical example of the oriented torus $\Sigma_{1,0}$ where the Teichm\"uller space $T(\Sigma_{1,0})$ is
defined to be the space of all flat metrics modulo similarity maps isotopic
to the identity. In his doctor thesis in 1913, J. Nielsen proved that
two homologous homeomorphism (resp. simple loops) of the torus are isotopic.
Thus the space of simple loops $\Cal S(\Sigma_{1,0})$ can be identified naturally with the set
of primitive elements in $H_1(\Sigma_{1,0}, \bold Z)$ modulo $\pm 1$ 
and the mapping class group $\Gamma (\Sigma_{1,0})$ is naturally isomorphic to
the automorphism group
$Aut^+(H_1(\Sigma_{1,0}, \bold Z))$. Define a \it marking \rm on $\Sigma_{1,0}$
to be a pair of oriented simple loops $(a, b)$ intersecting transversely
at one point. Fix a marking $(a,b)$ on $\Sigma_{1,0}$. Their homology  classes
$[a], [b]$ form a basis for the first homology group $H_1(\Sigma_{1,0},
\bold Z)$. In terms of the basis, one can identify $\Cal S(\Sigma_{1,0})$ with the set of
rational numbers $\hat \bold Q = \bold Q \cup \{\infty\}$ be sending
each primitive class $\pm (p[a]+q[b])$ to its ``slope" $p/q$. One can
also identify  the automorphism group $Aut^+(H_1(\Sigma_{1,0}, \bold Z))$ with $SL(2, \bold Z)$.
The natural action of the mapping class group $\Gamma (\Sigma_{1,0})$ 
on $\Cal S(\Sigma_{1,0})$ becomes the standard action of
$SL(2, \bold Z)$ on the rationals by fractional linear transformations.
The marking $(a,b)$ can also be used to parametrize the Teichm\"uller
space $T(\Sigma_{1,0})$ as follows. Fix a flat metric $d$ on $\Sigma_{1,0}$. We isotop $a$ and $b$ into two
$d$-geodesics $\hat a$ and $\hat b$ which intersect at one point $p$. 
Let $\theta$ be the angle measured from $\hat a$ to $\hat b$ at $p$
in the orientation of the surface and let $l_a$ and $l_b$ be the
lengths  of the geodesics $\hat a$ and $\hat b$. Assign to the flat
metric $d$ the complex number $z_d = l_a/ l_b e^{i \theta}$ in the
upper-half plane  $\bold H$. Evidently the invariant  $z_d$ depends
only on the similarity class of the flat metric $d$. Thus one obtains
a maps $\pi_m$ from the Teichm\"uller space $T(\Sigma_{1,1})$ to $\bold H$. This map
is a bijection since the inverse can  be constructed by sending $z \in
\bold H$ to the torus $\bold C /(\bold Z + z \bold Z)$ with marking
corresponding to $1$ and $z$. Note that the invariant $z_d$ is
independent of the orientations on $a$ and $b$. Furthermore,  the pair
$(z_d, z_d)$ is the Thurston  invariant of the parallelogram obtained by cutting
the flat torus $(\Sigma_{1,0}, d)$ open along the geodesics $\hat a$, $\hat b$.
Now if we are given a different marking $m' =(a', b')$, 
there is an $SL(2, \bold Z)$
matrix $A$ which sends $[a]$ to $[a']$ and $[b]$ to $\pm [b']$. A simple
calculation shows that two invariants $\pi_m$ and $\pi_{m'}$
are related by $A$ acting as the fractional linear transformation  on $\bold H$.
Thus  the Teichm\"uller space  $T(\Sigma_{1,0})$
  can be naturally identified with
$\bold H$ so that the action of the mapping class group becomes the
standard action of $SL(2, \bold Z)$ on $\bold H$ by fractional linear
transformations. In short, the three themes  $T(\Sigma_{1,0})$, $\Gamma (\Sigma_{1,0})$, and $\Cal S(\Sigma_{1,0})$
for the oriented  torus are exactly $(\bold H, SL(2,\bold Z), \hat \bold Q)$.
It is interesting to note that the complex structure on the Teichm\"uller
space makes both the angle function and the
 logarithm of ratio of length functions
pluriharmonic. Indeed, by fixing a marking $m =(a,b)$ on the torus 
$\Sigma_{1,0}$,
one obtains a fundamental domain map $f_m:  T(\Sigma_{1,0}) \to
T(4)/\bold R_+$ by sending the similarity class $[d]$ to the
parallelogram based on $\hat a$ and $\hat b$ which form a  fundamental
domain for the flat metric. The complex structure on $T(\Sigma_{1,0})$ makes
the map $f_m$ holomorphic, i.e., holomorphic motions in the Teichm\"uller
space corresponds to the homomorphic motions of the fundamental domains.
The same phenomenon does not seem to hold for the complex structure on the
Teichm\"uller space of a closed surface of higher genus with
hyperbolic metrics ([Ke1], [Wo2]). 
(A natural choice of fundamental domains for
closed surfaces of higher genus seems to be those associated to chains
in the surface. See Maskit [Ma] for more information.) However,
there are evidences indicating that 
the complex structure on the Teichmuller space is more closely
related to the singular flat metrics on the surface. 
To be more precise,
given a closed Riemann surface $\Sigma$ of genus $g$ and a fixed point on
the surface. There exists a unique singular
flat metric of area 1 in the conformal class of $\Sigma$ so that 
its sigularity is at the fixed point having cone angle $2\pi (2g-1)$.
See [Lu7] for more details on the relationship between these singular flat
metrics and the complex structure.

2.3. There is  a natural compactification of the upper-half plane $\bold H$ by
the real line $\hat \bold R = \bold R \cup \{\infty\}$ where the action
of $SL(2,\bold Z)$ extends continuously. Thurston's deep work on surface
theory shows that the same compactification also exists for all surfaces.
We shall discuss briefly Thurston's work in this section.
See [FLP] and [Th2] for more details.
To begin with, a proper 1-dimensional submanifold $s$ in a compact
surface $\Sigma$ is  called \it a curve  system \rm if no component of  $s$ 
is homotopic into $\partial \Sigma$ relative to $\partial \Sigma$. The
isotopy class of all curve systems on $\Sigma$ is denoted by $CS(\Sigma)$
and  was introduced by Dehn who called it \it the Arithmetic field \rm of the surface. In the
case of a torus, the set $CS(\Sigma)$ is naturally the set of all non-zero
lattice points in $H_1(\Sigma_{1,0}, \bold Z)$ modulo $\pm 1$. There
exists a quadratic pairing on $CS(\Sigma)$ given by the geometric intersection
number $I(\alpha, \beta) =$ min$\{|a \cap b|: a \in \alpha, b \in \beta\}$.
For the torus, the pairing is $I((p,q), (p',q')) = |pq' - p'q|$ which is the
absolute value of the canonical symplectic form on $\bold Z^2$. This
pairing satisfies the homogeneity and non-degenerate property in the
sense that $I(k_1 \alpha_1, k_2 \alpha_2) = k_1k_2 I(\alpha_1, \alpha_2)$
($k_i \in \bold Z_+$ and $k_i \alpha_i$ means $k_i$ copies of the curve system $\alpha_i$), and for each $\alpha$ there exists $\beta$ so that  $I(\alpha,
\beta) \neq 0$. Thurston's space of measured laminations $ML(\Sigma)$ is the
completion of $CS(\Sigma)$ with respect to the pairing $I$. In linear algebra,
given a non-degenerate quadratic form $\omega$ on a lattice 
$L$ of rank $r$, one can form a completion of $(L , \omega)$ by 
canonically embedding $L$ into $R^r$
so that the form $w$ extends continuously on $R^r$.
If the form is definite, the simplest way to construct the completion
is by formally extending $\omega$ to $\bold Q  L$ and taking the metric
completion of $\bold Q L$. If the form $\omega$ is not definite, one
may embed $L$ into the infinite dimensional space $\bold R^L$ (with the
product topology) by sending $x \in L$ to the linear function 
$\pi(x) = \omega(  . \quad, x)$. The canonical completion is given by taking the closure
of the set $\bold Q  \pi(L)$. Since the form $\omega$ is non-degenerate, the
Rieze representation theorem says that the closure is isomorphic to a vector
space $\bold R^r$ and the form $\omega$ extends continuously to the 
closure. Thurston's completion of $(CS(\Sigma), I)$ is the
analogous construction. The space $CS(\Sigma)$ is embedded into
 $\bold R^{\Cal S(\Sigma)}$
by sending $\alpha$ to the \it intersection  function \rm $Th(\alpha)
=I( . \quad, \alpha)$ and  the closure
of $\bold Q_+ Th(CS(\Sigma))$ is defined to be the completion, the 
space of measured lamination $ML(\Sigma)$. 
Thurston proved a remarkable theorem
that the space $ML(\Sigma)$ is homeomorphic to a Euclidean space and the quadratic
pairing  extends continuously to $ML(\Sigma)$ ([Bo2], [FLP], [Re],
[Th1]). Since Thurston's completion is
canonically constructed, the mapping class group $\Gamma(\Sigma)$ acts continuously
on $ML(\Sigma)$. In the case of torus, $ML(\Sigma)$ is canonically $H_1(\Sigma_{1,0}, \bold R)
/\pm 1$ and the action of the the mapping class group $SL(2, \bold Z)$
is the standard action. The projectivized  space PML 
$=(\Cal ML(\Sigma) -0)/\bold
R_{>0}$ is Thurston's compactification of the Teichm\"uller space.

\S3. {\bf  The Restriction Maps and the Reconstruction Theorems}

3.1.  A compact subsurface $\Sigma'$ in $\Sigma$ is called \it essential \rm 
if no component of $\partial \Sigma'$ is null homotopic in $\Sigma$. If
$\Sigma'$ is essential with negative Euler number, there exists a natural
restriction map from the Teichm\"uller spaces $T(\Sigma)$  to $T(\Sigma')$
(resp. from $ML(\Sigma)$  to $ML(\Sigma')$). The
restriction map is defined as follows. Given a hyperbolic metric
$d$ on $\Sigma$, we isotop the  open surface $int(\Sigma')$ (the interior
of $\Sigma'$) to an open subsurface $\Sigma''$ so that $\Sigma''$ is
bounded by disjoint simple geodesics. The metric completion of $(\Sigma'',
d|_{\Sigma''})$ is a hyperbolic metric $d'$ on $\Sigma'$ with geodesic
boundary. The restriction map sends $[d]$ to $[d']$.
The restriction map for the measured laminations  is defined similarly
(see [Lu2]). The key step is to define the restriction map from the space
of curve systems $CS(\Sigma)$ to $CS(\Sigma')$. Given $\alpha \in CS(\Sigma)$,
choose a representative $a \in \alpha$ so that the  number of
components of $a \cap \Sigma'$ is minimal.
Then the restriction map sends $[a]$ to $[a |_{\Sigma'}]$.
The restriction  maps are natural in the sense that if we are given two
essential subsurfaces $\Sigma_1 \subset \Sigma_2 \subset \Sigma$, then
the composition of restrictions is the restriction.

To state the reconstruction theorems, we  say that two essential subsurfaces
$\Sigma_1$ and $\Sigma_2$  \it overlap essentially \rm  if there is
a non-trivial simple loop which is isotopic into both $\Sigma_1$ and 
$\Sigma_2$. If furthermore  both $\Sigma_1$ and
 $\Sigma_2$ are level-1 subsurfaces, then their possible
intersection surfaces are either essential annuli, or an essential
 3-holed sphere or they are isotopic.

3.2.
With these preparation, we can state the main theorems in [Lu1], [Lu2]
and [Lu4] as follows which establish Grothendieck's reconstruction 
principle for the Teichm\"uller spaces, the measured  lamination spaces and
the mapping class groups.

{\bf Theorem}(Reconstruction of the Teichm\"uller spaces and the
measured lamination spaces). \it Each hyperbolic metric (resp. measured
lamination) on a surface of level at least 2 is constructed uniquely
up to isotopy by assigning a hyperbolic metric (resp. measured lamination)
to each essential level-1 subsurface so that when two level-1 subsurfaces
overlap essentially, the restrictions of the metrics (resp.
measured laminations) to their intersection are isotopic.
Furthermore, the restriction of the hyperbolic metric (resp. measured
lamination) on the surface to each level-1 essential subsurface
is isotopic to the assigned one. \rm

3.3.
For the mapping class group, it was a theorem of Dehn and Lickorish that
the mapping class group $\Gamma(\Sigma)$ is finitely generated by Dehn twists along
simple loops. The main result in [Lu4] states that

{\bf Theorem}(Reconstruction of the Mapping Class Groups). \it
Each orientation preserving self-homeomorphism of a surface which fixes the
boundary pointwise
 is isotopic
to a composition of finitely many Dehn twists so that  the composition
is unique modulo cancellation laws supported in subsurfaces of level-1.\rm

Earlier work on the subject were Hatcher and Thurston who  showed among
other things that the subsurfaces can be taken to be genus 2 with 3 holes
and Gervais who proved that subsurface can be taken to be genus 1 with
2 holes. Theorem 2 is just a simple simplification of the work of
Gervais.

3.4.
Before discussing the related theorems for surfaces of levels 0 or 1,
let us recall Thurston's embeddings of the Teichm\"uller space $T(\Sigma)$ and
the measured lamination spaces $ML(\Sigma)$. 
 Given an isotopy class $[d] \in$ $T(\Sigma)$, the
\it geodesic length function \rm $l_d : \Cal S(\Sigma) \to \bold R_{\geq 0}$ 
of the metric sends an isotopy class of simple loop to the length of
its geodesic representative.  For a measured lamination $m \in ML(\Sigma)$,
the \it  geometric intersection number function, \rm or simply
\it intersection function \rm  $I_m: \Cal S(\Sigma)
\to \bold R_{\geq 0}$ is given by $I_m(\alpha) = I(\alpha, m)$. Thurston's
embedding $Th: T(\Sigma) \to \bold R_{\geq 0}^{\Cal S(\Sigma)}$
(resp. $Th: ML(\Sigma) \to \bold R_{\geq 0}^{\Cal S(\Sigma)}$)
sends the isotopy class of a metric to its geodesic length function,
i.e., $Th([d]) = l_d$ 
(resp. sends a measured lamination to its intersection function).
The fact that the map $Th$ is injective for Teichm\"uller space was
a result of Fricke and Klein. The work of Okumura  [Ok1], [Ok2] 
and Schmultz [Sc]
determine the smallest finite set $F \subset \Cal S
(\Sigma)$ so that the restriction $l_d|_F$ determines the metric $d$.
The analogous question for measured lamination space seems to be
open. A result of Thurston shows $9g-9$ simple loops suffices to
determine an intersection function for closed
surface of genus $g$ ([Th1], [FLP]). But the number $9g-9$  is not the smallest.

3.5.
For the level-0 surface, i.e., the 3-holed sphere, the space of  simple
loops $\Cal S(\Sigma_{0,3})$ consists of isotopy classes of the three boundary components.
The Teichm\"uller space $T(\Sigma_{0,3})$, the measured lamination space $ML(\Sigma_{0,3})$ and
the mapping class group $\Gamma (\Sigma_{0,3})$ can be described as follows. By a theorem
of Fricke and Klein mentioned before, each isotopy class of a hyperbolic
metric on $\Sigma_{0,3}$ is determined by the lengths of three boundary components
and these lengths subject to no constraints, i.e., $Th(T(\Sigma_{0,3}))
= \bold R_{>0}^{\Cal S(\Sigma_{0,3})}$. The work of Thurston shows
that each measured lamination is determined by its intersection number
with three boundary components and these three numbers subject to no
constraints, i.e., $Th(ML(\Sigma_{0,3})) = \bold R_{\geq 0}^{\Cal S(\Sigma_{0,3})}$.  Dehn proved in 1938 that $\Gamma (\Sigma_{0,3}) $ is isomorphic to
free abelian groups on three generators which are the Dehn twists on
three boundary components.

With these results on 3-holed sphere, one can restate theorem  3.2 in an
equivalent form as follows. Given a  surface $\Sigma$ of level at least 2, a
real valued function on the set of simple loops $\Cal S(\Sigma)$ is a geodesic
length function (resp. an intersection function) if
and only if for each essential surface $\Sigma'$ of level-1, 
the restriction of the function to $\Cal S(\Sigma')$ is.

3.6.
For level-1 surfaces $\Sigma$, i.e., the 1-holed torus $\Sigma_{1,1}$ and 4-holed sphere $\Sigma_{0,4}$,
the set of simple loops $\Cal S(\Sigma)$, the Teichm\"uller space $T(\Sigma)$ and the
mapping class group $\Gamma(\Sigma)$ are essentially the  same as those of the torus. To be
more precise, let us consider the subset $\Cal S'(\Sigma)$ $\subset \Cal 
S(\Sigma)$ of isotopy
classes of simple loops which are not homotopic into the boundary
$\partial \Sigma$, i.e., $\Cal S'(\Sigma) = \Cal S(\Sigma) \cap CS(\Sigma)$.
There exists a natural bijection $i_*$ between $\Cal S'(\Sigma_{1,1})$ and $\Cal S(\Sigma_{1,0})$ induced
by the inclusion map from $\Sigma_{1,1}$ to $\Sigma_{1,0}$. This isomorphism preserves the
intersection pairing. For the 4-holed sphere $\Sigma_{0,4}$, there exists a
natural isomorphism $P^*$ between $\Cal S'(\Sigma_{0,4})$ 
and $\Cal S'(\Sigma_{1,1})$ which satisfies
$I(\alpha, \beta) = 2 I(P^*(\alpha), P^*(\beta))$. It is defined
as follows. Let $\tau$ be a hyperelliptic involution on the 1-holed
torus $\Sigma_{1,1}$ and let $P : \Sigma_{1,1} \to \Sigma_{1,1}/\tau$ be the
quotient map where $\Sigma_{1,1,}/\tau$ is the disc with three cone
points of order two (an orbifold). It is well known that the hyperelliptic
involution $\tau$ preserves the isotopy class of each simple loop and
$\tau$ commutes with each homeomorphism up to isotopy. Let the 4-holed
sphere $\Sigma_{0,4}$ be the subsurface of $\Sigma_{1,1}/\tau$ with three small
disc neighborhoods of the cone points  removed. Then the isomorphism
$P^*$ from $\Cal S'(\Sigma_{0,4})$ to $\Cal S'(\Sigma_{1,1})$ sends the isotopy class $[a]$ to $[b]$ where
$b$ is a component of $P^{-1}(a)$. To summarize, for a level-1 surface
$\Sigma$, there exists a bijection $\pi$ from $\Cal S'(\Sigma)$ to $\hat \bold Q$ so that
$\pi(\alpha) = p/q$ and $\pi(\beta) = p'/q'$ satisfy $pq'-p'q = \pm 1$
if and only if $I(\alpha, \beta) = 1$ for $\Sigma_{1,1}$ and 2 for $\Sigma_{0,4}$.
Draw a hyperbolic geodesic in the upper-half plane
ending at $p/q$ and $p'/q'$ when  $pq' -p'q = \pm 1$.
One obtains the so called ``modular configuration" (see figure 7).
Call three elements in $\Cal S'(\Sigma)$ forming an triangle 
if they correspond
to the vertices of an ideal triangle in the  modular configuration and call
four elements in $\Cal S'(\Sigma)$ forming a quadrilateral if
they correspond to the vertices of an ideal quadrilateral.
The modular structure on the space
of simple loops $\Cal S'(\Sigma)$ for level-1 surfaces was known
to Fricke and Klein ([Ke1]) and to Dehn ([De]) who used the rational
numbers to code the set $\Cal S'(\Sigma)$.
See also [HT], [MM],  [Se], [Th2] and others.

\midspace{0.1cm}
\centerline{\epsfbox{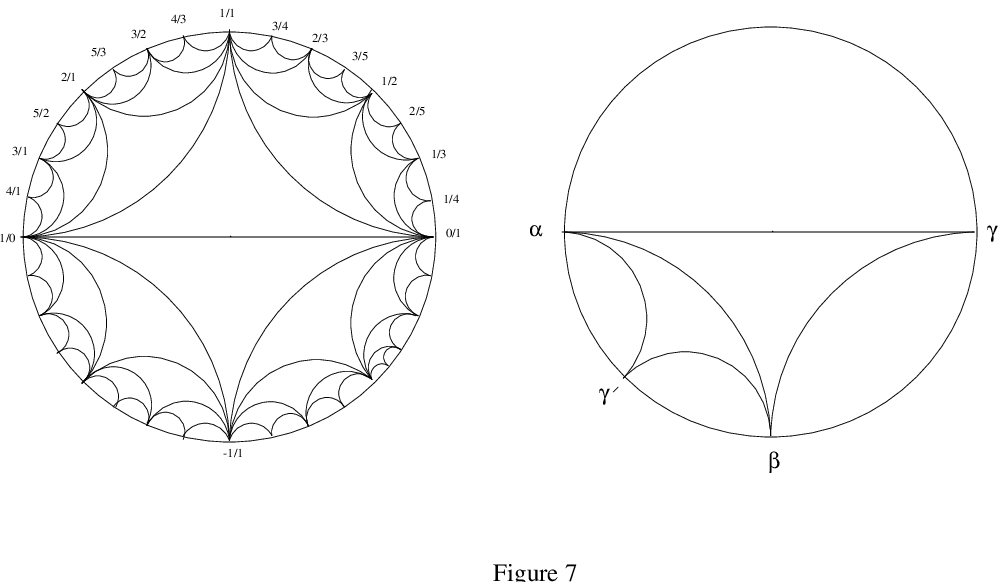}}

A special feature of the modular configuration
is the huge symmetry built in the configuration. This is, in our view,
one reason why the set of homotopy class of simple loops on the surface
is more useful than the fundamental group in establishing
the reconstruction principle for many structures.
Suppose we take four vertices $(\alpha, \beta, \gamma; \gamma')$
forming a quadrilateral in $\Cal S'(\Sigma)$ so that both
$(\alpha, \beta, \gamma)$ and $(\alpha, \beta, \gamma')$ are
triangles. Then there is an orientation reversing involution of 
$\bold QP^1$ leaving the quadrilateral fixed and interchanging
$\gamma$ and $\gamma'$. This involution is realized by an
orientation reversing involution of the surface $\Sigma$ which
is the reflection of the figure 9 (where $\gamma = \alpha \beta$) 
about the $yz$-plane.  On the other hand, given any triangle
$(\alpha, \beta, \gamma)$ in the modular configuration $\Cal
S'(\Sigma)$, there is a  $\bold Z_3$ action on  $\bold QP^1$
which permutes the three vertices. Thus there is a $\bold Z_3$
action on the surface $\Sigma$ permuting the isotopy classes.
This symmetry is illustrated in figure 8 below where the
1-holed torus is the Seifert surface of the trefoil knot and
the 4-holed sphere is the truncated boundary surface of a cube.
The symmetry involved in the 4-holed sphere is huge which is
difficult to visualize in figure 9. Indeed, as one can see from
figure 8, any permutation of the four boundary components is
realized by a homeomorphism preserving the set \{$\alpha,
\beta, \gamma\}$, i.e., the permutation group on four letters
acts on $\Sigma_{0,4}$ preserving  set \{$\alpha,
\beta, \gamma\}$. 

\midspace{0.1cm}
\centerline{\epsfbox{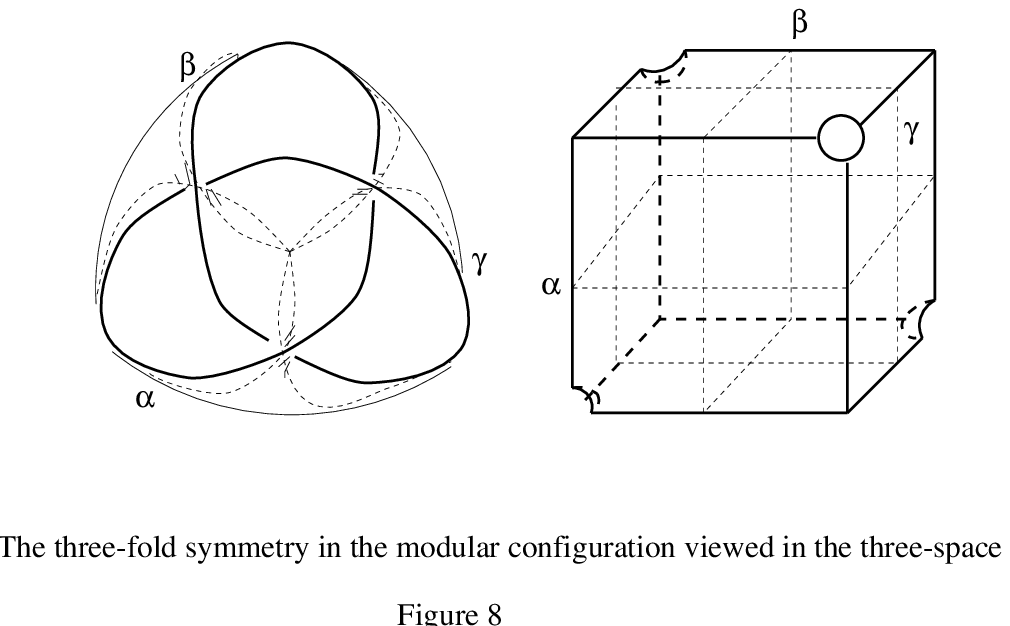}}

As an application of this 24-fold symmetry, we consider the
trace relations for $SL(2, \bold C)$ matrices. In this case, the
analogous question to problem 2 in \S 1.1 for triangles and quadrilaterals
is the following. Given three complex numbers $a,b,c$, do there
exist two matrices $A, B \in SL(2, \bold C)$ so that $tr(A) = a,
tr(B) = b$ and $tr(AB) = c$? It is well known that the answer is
positive. Next, in analogous to six edge lengths of a quadrilateral,
given seven complex numbers $a_1, a_2, a_3, a_{12}, a_{23}, a_{31}$
and $a_{123}$, under what condition do there exist three
$SL(2, \bold C)$ matrices $A_1, A_2,$ and $A_3$ so that
$tr(A_i) = a_i$, $tr(A_iA_j) = a_{ij}$ and $tr(A_1A_2A_3) = a_{123}$?
A solution by Fricke-Klein [FK] and Vogt [Vo] was the following. 
These three matrices exist if and only if
$a_{123}^2 - a_{123}(a_1 a_{23} + a_{2}a_{31} + a_3 a_{12} -a_1a_2a_3)
a_1^2 + a_2^2 a_3^2 + a_{12}^2 + a_{23}^2 + a_{31}^2 + $$
a_{12}a_{23}a_{31}
+a_1a_2a_3 + a_{12}a_{23}a_{31} - a_1a_2a_{12} - a_2a_3a_{23} - a_3a_1a_{31}
-4 =0$. This equation, as it stands, is quite complicated. One can
easily notice the 3-fold symmetry of the equation under cyclic permutation
of $\{a_1, a_2, a_3\}$. In fact, there exists 24-fold symmetry in the
equation. Namely, the equation is  invariant under any permutation
of $\{a_1, a_2, a_3, a_{123}$\}. This can be seen using the 
modular configuration. Indeed, if we choose the generators of the 
fundamental group of the 4-holed sphere carefully (see for instance, 
figure 5 in [Lu1]), then 
the four boundary components are represented by $x_1, x_2, x_3, 
x_1x_2x_3$ 
and three simple loops forming a triangle in the modular configuration
by $x_1x_2$, $x_2x_3$ and $x_3x_1$.  The first equation in Theorem A4 (b) 
in the appendix A is a rewriting of the above polynomial  equation
in terms of the modular configuration.

3.7. The Teichm\"uller space, measured lamination space and the
mapping class groups for level-1 surfaces can be explicitly constructed
from the modular configuration on the set of simple loops.
To be more precise, the geodesic length functions, the geometric intersection
number functions and the Dehn twists satisfy universal relations at the
vertices of triangles and quadrilaterals in $\hat \bold Q$. 
Furthermore, these universal
relations form a complete set of relations. See Appendix A for the
list of universal relations.

The relations for the Dehn twists were found by Dehn  in [De]. D. Johnson [Jo]
independently rediscovered the lantern relation (relation (IV) in Theorem
A3 in the Appendix A) in 1979.  Dehn also
proved that these relations are complete for the mapping class
group of level-1 surfaces. The relations for the geodesic length function
were essentially discovered by Fricke and Klein [FK] and Vogt [Vo]
(thought they were not stated in terms of the modular relations). These are derived
from the trace identities for SL(2, $\bold C$) matrices. That the set  of
all relations is complete was proved by Keen [Ke2] for the 1-holed torus
and was proved in [Lu1] for 4-holed spheres using Maskit combination theorem.
The relations for the measured laminations are just the degenerations
of that of hyperbolic metrics and they are shown to be complete in [Lu2].

3.8.
The relationship between the Teichm\"uller spaces and the mapping class
groups among the level-1 surfaces becomes much clear if one considers
the Teichm\"uller spaces $T_{1,1}$ and $T_{0,4}$
of complete hyperbolic metrics
with cups ends (on the open surface), and the reduced mapping class group
$\Gamma^*(\Sigma)$ which is the quotient of the mapping class group
by the subgroup generated by Dehn twists on boundary components.
The key fact is that the hyperelliptic involution $\tau$ on
$\Sigma_{1,1}$ induces the identity map on both the Teichm\"uller space
$T_{1,1}$ and $\Cal S(\Sigma_{1,1})$ and is in the center of the
mapping class group.
Indeed, one has a  natural  biholomorphism
between $T_{1,1}$ and $T_{0,4}$ induced by the pull back map
$P: \Sigma_{1,1} \to \Sigma_{1,1}/\tau$.  A natural isomorphism
from $\Gamma^*(\Sigma_{1,1})$ to $\Gamma(\Sigma_{1,0}) = SL(2, \bold Z)$
is
induced by inclusion of $\Sigma_{1,1}$ to $\Sigma_{1,0}$.
Since the hyperelliptic
involution $\tau$ commutes with each homeomorphism, there is an
monomorphism from the reduced mapping class group $\Gamma^*(\Sigma_{0,4})
\to \Gamma^*(\Sigma_{1,1})/ <\tau> = PSL(2, \bold Z)$ whose image is the
principal congruence subgroup of order 2.

3.9.
It is instructive to read [Gr1] on the related topics. We cite the
paragraph on page 248-249 in [Gr1] below (with English translation by
L. Schneps).
`` There is a striking analogy, and I am certain it is not merely
formal, between this principle and the analogous
principle of Demazure for the structure of reductive algebraic groups,
if we replace the term `level' or `modular dimension' with 
`semi-simple rank of the reductive group'. The link becomes even
more striking, if we recall  that the mapping class group  $\Gamma^*_{1,1}$
is no other than $SL(2, \bold Z)$, i.e., the  group of integral points
of the simple group scheme of `absolute' rank 1 $SL(2)_{\bold Z}$.
Thus, $\underline{\text{the fundamental building block for
the Teichm\"uller tower is essentially the same as for the}}$ \newline
 $\underline{\text{ 
tower of reductive groups of all ranks}}$ - a 
group of which, moreover, we may say that it is doubtless
in all the essential disciplines of mathematics."


\S 4. {\bf The Space of Simple Loops on Surfaces and the Modular Structure}

4.1. Unlike subsurfaces in a surface, simple loops on surfaces have been
the focus of more attention for a long time. Indeed, most of the surface
problems can be reduced to ones concerning simple loops and the
proofs of theorems 3.2 and 3.3 are no exception. The topological investigation
of the set $\Cal S(\Sigma)$ of isotopy classes of essential simple loops began in Dehn's
work [De] on the mapping class groups. As an example of use of simple loops
to solve surface problems, let us recall the elegant proof of Dehn that
the (reduced) mapping class group $\Gamma^*(\Sigma_{0,4})$ is a free group on two generators
generated by Dehn twists on two simple loops intersecting at two points.
Dehn first observed that the set $\Cal S'(\Sigma_{0,4})$ of essential simple loops not homotopic
into the boundary forms the modular configuration $\hat \bold Q$ and the mapping
class group $\Gamma^*(\Sigma_{0,4})$ acts on the modular configuration faithfully preserving
both the modular relation and the orientation. Thus $\Gamma^*(\Sigma_{0,4})$ is a subgroup of
the modular group $PSL(2, \bold Z)$. Since each boundary component of 
$\Sigma_{0,4}$ is fixed by the mapping class group elements, $\Gamma^*(\Sigma_{0,4})$ is actually in the
principal congruence subgroup of level 2 generated by two matrices
$\left(\smallmatrix 1&2\\0&1\endsmallmatrix\right)$ and
$\left(\smallmatrix 1&0\\-2&1\endsmallmatrix\right)$.
But these two matrices correspond to two Dehn twists mentioned above.

To go from simple loops to subsurfaces, one takes the regular neighborhood
of a union of simple loops. In this way, it can be shown for instance that
given any two level-1 essential subsurfaces $A, B$, there  is a sequence of level-1
essential subsurfaces starting from $A$ and ending at $B$ so that any two
adjacent level-1 subsurfaces overlap in an essential level-0 subsurface.

4.2. The work of Dehn [De] and Lickorish [Li] already suggested strongly
that level-1 subsurfaces are fundamental in simplifying the intersections
of two simple loops. Indeed, lemma 2 in [Li] states that if two simple
loops $a, b$ satisfy either $|a \cap b| \geq 3$ or $|a \cap b| = 2$ with
non-zero algebraic intersection number, then there is a Dehn twist which
sends $b$ to a new loop having fewer intersection points with $a$. Thus
the only situation which cannot be simplified are: 1) $a, b$
are disjoint, 2) $a$ intersects $b$ at one point and 3) $a$ intersects
$b$ at two points of different intersection signs.  In these cases, the
lowest level connected subsurface which contains both $a$ and $b$  is 
either a level-0 or a level-1 subsurface. Furthermore, the pair of curves
$a, b$ satisfying conditions 2) or 3) corresponds to the basic relation
in the modular configuration.  

4.3.  We have mentioned in several places the notion of modular structure.
Here is a formal definition after Thurston's geometric structures on manifolds.

\noindent
{\bf Definition.} (a)
 A ($\hat \bold Q$, $PSL(2, \bold Z))$ \it  modular structure \rm on
a set $X$ is a maximal collection of charts $\{(U_i, \phi_i) | \phi_i:
U_i \to \hat \bold Q$ is injective$\}$ so that the following hold.

(1)  $X = \cup_i U_i$.

(2) The transition function $\phi_i \phi_j^{-1}$ is the restriction of
an element in $PSL(2, \bold Z)$.

(b) A modular structure on a set $X$ is called \it compact \rm
if the group of bijections of X which preserve the modular structure
acts on $X$ with finite orbits.

It seems that compactness is essential for developing a useful
``function theory" on a set with a modular structure.
All interesting examples that we encounter have compact modular structures.

For an oriented surface $\Sigma$ of level at least 1, the set $\Cal S'(\Sigma)$ of isotopy
classes of non-boundary parallel essential simple loops on $\Sigma$ has a 
natural compact
modular structure invariant under the action of the mapping class group.
A special collection of  charts for the modular structure is 
given by $(\Cal S'(\Sigma'),
\phi_{\Sigma'})$ where  $\Sigma'$ is an essential level-1 subsurface
and $\phi_{\Sigma'}: \Cal S'(\Sigma')
\to \Cal S'(\Sigma_{1,1}) = \hat \bold Q$ is a bijection induced by
either an orientation preserving homeomorphism or by an orientation 
preserving quotient map (see \S3.6). To see that the second condition (2)
in the definition holds, one simply notes that if two essential level-1 subsurfaces
intersect at two non-homotopic simple loops, then they are isotopic.
To see the compactness, we note that the mapping class group acts
on $\Cal S'(\Sigma)$ preserving the modular structure and the
action of the mapping class group has finite orbits. Thus
we can
talk about triangle and quadrilateral in $\Cal S'(\Sigma)$. Furthermore, since 
the set of rational numbers $\hat \bold Q$ has the natural 
orientation invariant under $PSL(2, \bold Z)$, we 
can talk about oriented triangles in $\Cal S'(\Sigma)$. 

Another example of compact modular structure 
is the set of isotopy classes of 3-holed
sphere decompositions of a surface which seems to be related to the Heegaard
splittings of 3-manifolds. See Appendix B for more detail.

4.4. One way to see the modular structure on the space of simple loops $\Cal S(\Sigma)$
is to use the notion of resolution of intersection points. Recall that
two rational numbers $p/q$ and $p'/q'$ are modular related if $pq'-p'q =
\pm 1$ and are denoted by $p/q \perp p'/q'$. Two isotopy classes
$\alpha$ and $\beta$ of curves
on surfaces corresponding to a  modular related pair
are denoted by $\alpha \perp \beta$ or  $\alpha \perp_0 \beta$.
Here $\alpha \perp \beta$ means that their intersection number
 $I(\alpha, \beta) =1$ and $\alpha \perp_0
\beta$ means that  $I(\alpha, \beta)= 0$ with zero algebraic intersection number.  
To find out the vertices of ideal triangles based on $\{\alpha, \beta\}$,
we use the resolutions of intersections. Recall that surfaces are oriented.
If $a, b$ are two arcs intersecting at one point transversely, then the
\it resolution of $a \cup b$ at the intersection point from $a$ to $b$  \rm
is
defined as follows. Fix any orientation on $a$ and use the orientation on
the surface to determine an orientation on $b$. Then resolve the intersection
according to the orientations (see figure 8). The resolution is independent
of the orientation chosen on $a$. If $\alpha \perp \beta$ or
$\alpha \perp_0 \beta$, take $a \in \alpha$ and $b \in \beta$ so that
$|a \cap b| = I(\alpha, \beta)$. Then the curve obtained by resolving
all intersection points in $a \cap b$ from $a$ to $b$ is again an essential
non-boundary parallel simple loop. 
We denote the isotopy class by $\alpha \beta$.
One sees easily that positively oriented triangles  and  quadrilaterals
in the modular structure on $\Cal S'(\Sigma)$ are exactly $(\alpha, \beta, \alpha\beta)$
and $(\alpha, \beta, \alpha \beta, \beta \alpha)$. If $\alpha \perp \beta$
or $\alpha \perp_0 \beta$, we use $\partial(\alpha \cup \beta)$ to denote
the isotopy class of the boundary of a regular neighborhood of $a \cup b$.
In terms of these notations, all universal relations for the geodesic
length functions, the intersection functions and the Dehn twists are
expressed in terms of $\alpha, \beta, \alpha\beta, \beta\alpha$ and
the components of $\partial (\alpha \cup \beta)$. For instance, the
relations for the Dehn twists are: 1) if $\alpha \perp \beta$, then
$D_{\alpha} D_{\beta} = D_{\beta} D_{\alpha \beta}$ and $(D_{\alpha}
D_{\beta} D_{\alpha \beta})^4 = D_{\partial (\alpha \cup \beta)}$,
and 2) if $\alpha \perp_0 \beta$, then  $D_{\alpha}D_{\beta}D_{\alpha \beta}
= D_{\partial(\alpha \cup \beta)}$ (the lantern relation).
Since the modular relation $(\hat \bold Q, \perp)$ has a $\bold 
Z_3$-symmetry leaving an ideal triangle invariant, we obtain
$\alpha (\beta \alpha) = (\alpha \beta) \alpha = \beta$.

\midspace{0.1cm}
\centerline{\epsfbox{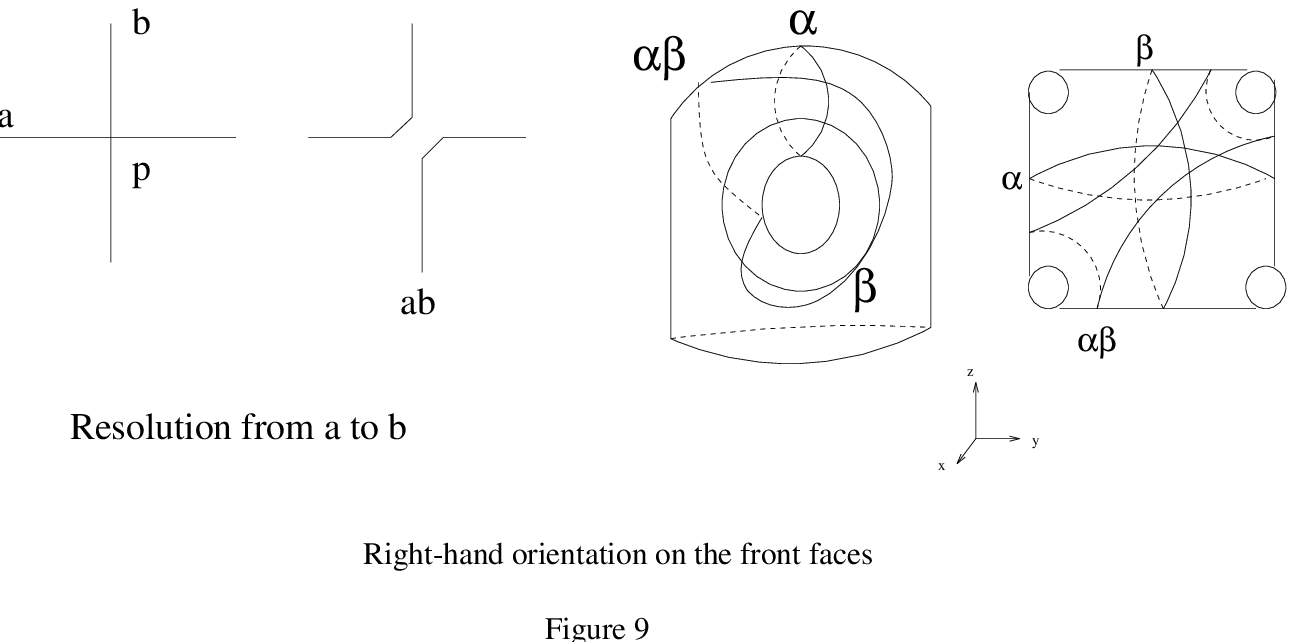}}
\midspace{0.1cm}

4.5. One of the most useful property of the modular structure  on $\Cal S'(\Sigma)$ is the following
lemma (lemma 7 in [Lu1]) which generalizes Lickorish's lemma 2 in [Li].
It states that given two intersecting elements $\alpha, \beta \in $ $\Cal S'(\Sigma)$ which are
not related by the modular relation $\perp$ or $\perp_0$, then we can write
$\beta = \gamma_1 \gamma_2$ with $\gamma_1 \perp \gamma_2$ or $\gamma_1 \perp_0
\gamma_2$ so that (1) $I(\alpha, \gamma_i) < I(\alpha, \beta)$ and
$I(\alpha, \gamma_2 \gamma_1) < I(\alpha, \beta)$ for $i=1,2$ and
(2) if $\gamma_1 \perp_0 \gamma_2$, then for each component $\delta$ of
$\partial (\gamma_1 \cup \gamma_2)$ we have $I(\alpha, \delta) < I(\alpha, \beta)$.
As an easily consequence, one shows that the reconstruction principle
for the Teichmuller spaces follows from theorem 3.2 for level-2 surfaces.
As another consequence, one shows that the space $\Cal S(\Sigma)$ is finitely generated
in the following strong sense: There is a finite  subset $X_0$ in $\Cal S(\Sigma)$
so that $\Cal S(\Sigma)$ = $\cup_{n=0}^{\infty} X_n$ where $X_{i+1} = X_i
\cup \{ \alpha | \alpha = \gamma_1 \gamma_2 $ where either (1) $\gamma_1 \perp
 \gamma_2$, and $\gamma_1, \gamma_2, \gamma_2\gamma_1$ are in $X_i$ or
(2) $\gamma_1 \perp_0 \gamma_2$ and $\gamma_1, \gamma_2, \gamma_2 \gamma_1$
and each component of $\partial (\gamma_1 \cup \gamma_2)$ are in $X_i$\}.

\midspace{0.1cm}
\centerline{\epsfbox{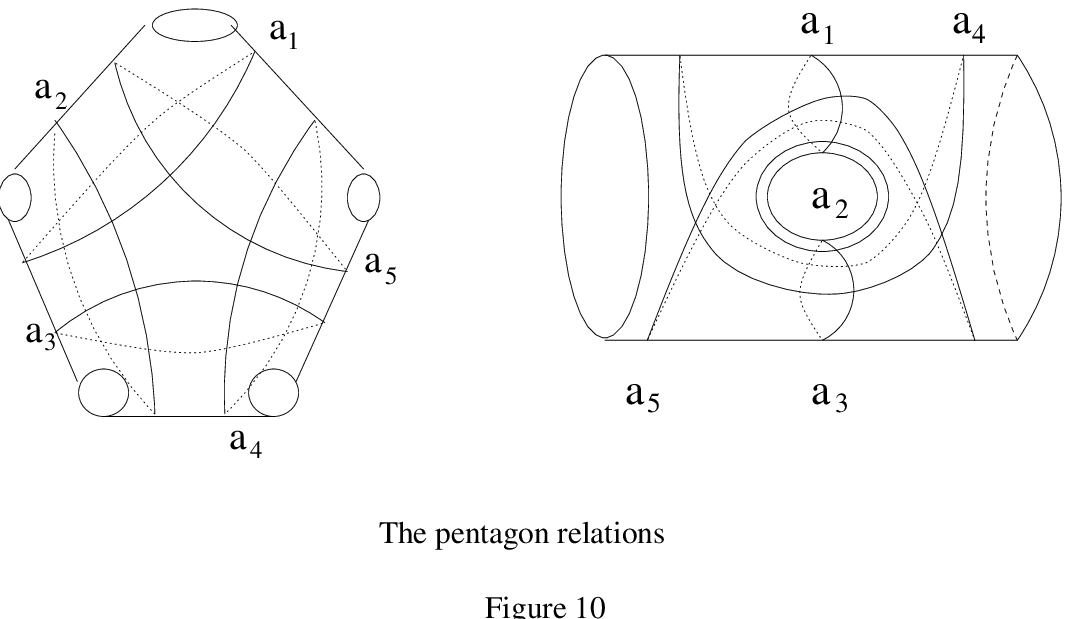}}

4.6. The proof of the reconstruction theorem for level-2 surfaces has always
been one of the key steps in establishing the reconstruction principle.
In dealing with
simple loops on leve-2 surfaces, the following collection of five curves
(the pentagon relation) $\{\alpha_1, ..., \alpha_5 $ in  $\Cal S'(\Sigma) | I(\alpha_i, \alpha_j)
= 0 $ for indices $|i -j| \neq 1$ mod(5)$\}$ (see figure 10) appears constantly
and plays an important role. For 5-holed sphere, one has the following
relation $\alpha_i \alpha_{i+1} \alpha_{i+2} = \alpha_{i+3} \alpha_{i+4}$
and for the 2-holed torus, we have $\alpha_1 \alpha_2 \alpha_3 \alpha_4
=\alpha_3 \alpha_2 \alpha_1$ (see [Lu2]). These five curves for the 5-holed
sphere were  first observed by Dehn in [De] who showed that the Dehn twists
on them generate the reduced mapping  class group for  5-holed sphere.
Furthermore,  these five curves are rigid in the sense
that any other collection of five curves with the same disjointness property
is the image of $\{\alpha_i\}$ under a homeomorphism ([Lu5]).

4.7. We finish this section with an application of the notion of resolving
intersection to a multiplicative structure on the space of curve systems
$CS(\Sigma)$. 
Given two curve systems $a, b$ on an oriented surface with
$|a \cap b| = I([a], [b])$, the multiplication $ab$ is defined to be the
1-dimensional submanifold obtained by resolving all intersection points
in $a \cap b$ from $a$ to $b$.  It can be shown that $ab$ is again a
curve system (see Appendix C for a simple proof when $a, b$ contain
no arcs).  This induces a multiplicative structure on $CS(\Sigma)$ by defining
$\alpha \beta = [ab]$ where $a \in \alpha$, $b \in \beta$ and
$|a \cap b| = I(\alpha, \beta)$.  For instance the Dehn twist on a
simple loop $\alpha$ applied to $\beta$ is given by $D_{\alpha}(\beta)
= \alpha^k \beta$ where $k =I(\alpha, \beta)$. This multiplication is
natural with respect to the action of the mapping class group and is highly non-commutative.
Indeed, if $\alpha$ contains no arc component, then $\alpha \beta = \beta \alpha$
implies $I(\alpha, \beta) =0$. As a consequence of this, one obtains a new
proof of a result of Ivanov ([Iv]) that Dehn twists on two intersecting isotopy
classes of  simple
loops can never be commuting up to  isotopy. The most interesting
property of the multiplication seems to be the ``cancellation law" that
if each component  of $\alpha$ is not an arc and intersects $\beta$, than
$\alpha(\beta\alpha) = (\alpha \beta) \alpha = \beta$. This is a generalization
of the $\bold Z_3$-symmetry in the modular configuration. As an application
of the cancellation law, let us prove a weak form of a result of Thurston
[Th2] that if $\alpha$ and $\beta$ are two surface
filling simple loops (i.e., $I(\alpha, \gamma) + I(\beta, \gamma) > 0
$ for all $\gamma \in$ $\Cal S'(\Sigma)$), then the self-homeomorphism $f = D^{-1}_{\alpha}
D_{\beta}$ does not leave any curve system invariant up to isotopy
(in fact Thurston proved that $f$ is pseudo-anosov). Indeed, if $f$ leaves
an element $\gamma \in$ $CS(\Sigma)$ invariant, then $D_{\alpha} (\gamma)
=D_{\beta}(\gamma)$, i.e., $ \alpha^k \gamma = \beta^l \gamma$. Now multiply
the equation by $\gamma$  from the left and use the cancellation law.
One obtains a contradiction to the surface filling property.

\S5. {\bf Reduction to Level-2 Surfaces}

The goal of this section is to establish a fairly general criterion
to reducing problems concerning all surfaces to that of level-2
surfaces. 

5.1.
We shall  begin by some abstract definitions. Given a  subset $X$
of the set $Y^{\Cal S(\Sigma)}$ of all maps from $\Cal S(\Sigma)$ to $Y$, 
we say that the subset $X$  has \it property RP \rm if for any decomposition
of  the surface $\Sigma$ as a union
of two essential subsurfaces $A_1$ and $A_2$ which overlap in a level-0
essential subsurface, then the restriction map from $Y^{\Cal S(\Sigma)}$ to
$Y^{\Cal S(A_1) \cup \Cal S(A_2)}$ is injective on $X$. In another  words,
if $f, g$ are two elements in $X$ so that $f|_{\Cal S(A_1) \cup \Cal S(A_2)}
= g|_{\Cal S(A_1) \cup \Cal S(A_2)}$, then $f = g$. For simplicity,
we call the function $f|_{\Cal  S(A_i)}$  the \it restriction \rm
of $f$ to $A_i$. We say a subset 
$X \subset Y^{\Cal S(\Sigma)}$ with property RP is \it complete \rm if
for any two elements $f_1$ and $f_2$ in the restriction of $X$ to $A_1$ and 
$A_2$ so that  their restrictions to the overlap  $A_1 \cap A_2$ are the same, 
then there exists an element $f \in X$ whose restrictions 
to $A_i$ is $f_i$ for $i=1,2$.
For instance, the set of all geodesic length functions and the
set of all intersection functions have  complete property RP. 
 This is equivalent to the following gluing lemma for
hyperbolic metrics and measured laminations. Namely, suppose
the surface $\Sigma$ is a union of
two essential subsurfaces $A_1$ and $A_2$ which overlap in an essential
level-0 surface. If  we are
given two hyperbolic metrics $d_i$ on $A_i$ whose restrictions to the
overlap  of $A_1$ with $A_2$ are isotopic, then there is a hyperbolic
metric unique up to  isotopy  on the surface $\Sigma$ 
whose restriction to $A_i$ is isotopic
to $d_i$. The same gluing lemma holds for measured laminations.
Note that SL(2, $\bold C)$ characters 
do not have property RP due to the existence of reducible representations.
 The mapping
class group $\Gamma(\Sigma)$ considered as a subset of $\Cal S(\Sigma)^{\Cal
S(\Sigma)}$  does not have  property RP either.  But if one modifies the
definition of $\Cal S(\Sigma)$ by taking the isotopy classes of
all \it oriented \rm   simple loops, then the mapping class group has
 complete property RP.

5.2.  The main reduction lemma says the following. If $X$ is a subset of
$Y^{\Cal S(\Sigma)}$ so that for each level-2 essential subsurface
$\Sigma'$ the restriction of $X$ to $Y^{\Cal S(\Sigma')}$ has property
RP, then $X$ has property RP. See Appendix D for a proof of this
reduction lemma.
As a consequence of this reduction lemma, we have the following
fact. Suppose $\Sigma$ is a surface of level at least three and
$X$ and $X'$ are two subsets of $Y^{\Cal S(\Sigma)}$ so that
for each level-2 essential subsurface $\Sigma'$ the restrictions
of $X$ and $X'$ to $Y^{\Cal S(\Sigma')}$ are the same.
If furthermore that $X \subset X'$ and $X$ has complete property RP,
then $X = X'$.
 To see this, we use induction on the level of subsurfaces.
First of all, by the reduction lemma, both $X$ and $X'$ have property
RP. Now to show $X' \subset X$, take an element $x' \in X'$ and
decompose $\Sigma$ into a union of two essential surfaces $A_1$ and
$A_2$ of smaller levels so that they overlap in a level-0 surface.
By the induction hypothesis, we find two elements $x_1$ and $x_2$
which are in the restrictions of  $X$ to $A_1$ and $A_2$ so that
$x_i$ is the restriction of $x'$ to $A_i$. But the restrictions
of $x_i$ to the overlap are the same, namely, it is the restriction of $x'$
to the overlap. Thus by the completeness, there is an  element $x$ in
$X$ whose restrictions to $A_i$ is $x_i$, Thus $x = x'$ by property RP.
This shows $X' \subset X$.

By taking  $X$ to be the set of all geodesic length
functions and $X'$ to be the set of all real valued functions on $\Cal 
S(\Sigma)$ which satisfy the universal relations in theorem A1 in Appendix
A, we see that the reconstruction principle for all Teichmuller spaces
follows from that for level-2 surfaces.

Also, the reduction lemma shows that the problem on
the  automorphisms of the curve complex of a surface is essentially
a problem on level-2 surfaces (see [Lu5]).

5.3. The above gives some hints on the special role played by 
level-2 surfaces.  It also supports 
Grothendick's principle that in the reconstruction process
``relations are supported in level-2 surfaces".

\S6. {\bf SL(2, $\bold C)$ Representation Variety of Surface Groups}

6.1. An $SL(2, \bold C)$ representation of a group is a homomorphism of
the group into $SL(2, \bold C)$. The character of the representation
sends each group element to the trace of the representation matrix. 
If the group is the fundamental group of a surface, by using
a result of Fricke and Klein [FK] and Vogt [Vo], one shows that
the character function is determined by its
restriction to the set $\Cal S(\Sigma)$ of homotopy  classes of simple loops.
The main result in [Lu3] shows that the character function on $\Cal S(\Sigma)$
satisfies the reconstruction principle, i.e., 
except for finitely many (at least $2^{n-1}$) exceptional functions
defined on $\Cal S(\Sigma_{0,n})$  for $n \geq 5$,
a function on $\Cal S(\Sigma)$ is
an $SL(2, \bold C)$ character if and only if for each 
essential level-1 subsurface $\Sigma'$ in $\Sigma$ the 
restriction of the function to  $\Cal S(\Sigma')$ is
an $SL(2, \bold C)$  character.
An exceptional function  $f: S(\Sigma_{0,n}) \to 
\bold C$ satisfies the
following (1) $f(S(\Sigma_{0,n})) = \{2, -2\}$, (2) for each level-1
subsurface, the restriction of $f$ to the subsurface is a character,
(3) there exists
a level-2 subsurface $\Sigma'$ so that $f|_{\Cal S(\Sigma')}$
is exceptional. 
All exceptional functions
are constructed from the basic one defined on
$S(\Sigma_{0,5})$ which sends $b_i$ to 2 and all others to $-2$.
There are no 
representations whose characters are these exceptional functions.

Given a surface 
$\Sigma$ of level-1, $SL(2, \bold C)$ characters
on $\Cal S(\Sigma)$ are characterized by the trace identities on 
the vertices of triangles
and quadrilaterals in the modular relation (see Appendix A for the exact
statement). Thus the space of all $SL(2, \bold C)$ characters of a surface group
can be explicitly described. The reconstruction theorem also holds
for any $SL(2, K)$ characters where the field $K$ is quadratically 
complete (i.e., each quadratic equation with coefficients in $K$
has roots in $K$).

6.2. The main difficulty in establishing the reconstruction principle
for $SL(2, \bold C)$ characters is due to the existence of reducible representations.
Recall that an $SL(2, \bold C)$ representation is \it reducible \rm if it leaves a
1-dimensional linear subspace in $\bold C^2$ invariant.  Unlike the
discrete faithful subgroups in $SL(2, \bold R)$ which occur in the
Teichmuller spaces, there are many irreducible representations of a
surface group so that its restriction to a subgroup coming from
an essential subsurface of negative Euler number is reducible. 
Now the gluing
lemma (see \S5.2) is valid only for representations so that their restrictions
to the fundamental group of the intersection surface are irreducible.
Thus one should choose the decomposition of a surface $\Sigma$ as a union
of two subsurfaces $\Sigma_1$ and $\Sigma_2$ carefully. It turns out 
that the following is true 
which plays a key role in chosing the decomposition of a surface. Namely,
a representation of a surface group into $SL(2, \bold C)$ is irreducible if and
only if its restriction to the subgroup of a Euler number -1 subsurface is
irreducible ([Lu3]).

As a consequence, we obtain the following result concerning 
$SL(2,K)$ characters on any group. Suppose $K$ is a
field so that each quadratic equation with coefficients in  $K$
has roots in $K$.  Given a  group $G$, 
we are interested in finding all $SL(2,K)$ characters on
$G$.  In his work on $SL(2,\bold R)$ characters, Helling [He]
introduced the  notion of \it trace function. \rm Recall that 
a $K$ valued function $f$ on $G$ is a \it
trace function \rm if
any  two elements $x, y$ in $G$, $f(xy) + f(xy^{-1}) = f(x)f(y)$ and
$f(id) = 2$. 
Evidently all $SL(2,K)$ characters on $G$ are trace functions
due to the trace identity $tr(AB) + tr(AB^{-1}) = tr(A) tr(B)$.
One consequence of the characterization theorem is that
each trace function is also a character.

6.3. The role of level-1 surfaces among all surfaces is similar to the 
role of 2-generator groups among all groups.  For instance, by Jorgensen's 
inequality, a non-elementary subgroup in $SL(2, \bold C)$ is discrete 
if and only if each of its 2-generator subgroup is discrete. The 
reconstruction theorem for the Teichmuller space says that a faithful
 representation of a surface group into $SL(2, \bold R)$ is discrete if 
and only if its restriction to each subgroup of its level-1 subsurface 
is discrete and uniformizing a surface of the same type. It is
natural to ask if the similar description exists for discrete close
surface subgroup of $SL(2, \bold C)$.

6.4. It is interesting to ask if the reconstruction principle holds for 
representations of the surfaces groups into  general linear group 
$GL(n, \bold C)$.  To be more precise, suppose $f$ is a complex valued 
function defined on the conjugacy classes of the fundamental
 group of the surface so that the restriction of $f$ to the conjugacy classes 
of the fundamental group of each essential level-1  subsurface is a 
$GL(n, \bold C)$-character. Is $f$  the character of some 
$GL(n, \bold C)$  representation of the surface group?

\noindent
{\bf Appendix A.} The Statement 
of the Reconstruction Theorems for Level-1 Surfaces 

\medskip

Given a hyperbolic metric $d$ on a surface, the \it trace \rm of
the metric $d$ is the function $2 cosh l_d/2$ where $l_d$ is the
geodesic length function associated to $d$.

\noindent
{\bf Theorem A1.} \it (a) For the surface $\Sigma_{1,1}$  with $b =\partial \Sigma
_{1,1}$,
a function $t: \Cal S(\Sigma_{1,1}) \to \bold R_{\geq 2}$  is a trace
function of a hyperbolic metric  if and only if
the following hold.
$$ \prod_{i=1}^3 t(\alpha_i) = \sum_{i=1}^3 t^2(\alpha_i) + t(b) -2 \quad
\text{and}$$
$$ t(\alpha_3) t(\alpha_3') =  \sum_{i=1}^2 t^2(\alpha_i) + t(b) -2$$
where $(\alpha_1, \alpha_2, \alpha_3)$ and $(\alpha_1, \alpha_2, \alpha_3')$
are distinct ideal triangles in $\Cal S' (\Sigma_{1,1})$.

(b) For the surface $\Sigma_{0,4}$ with $\partial \Sigma_{0,4} = \cup_{i=1}^4 b_i$,
a function 
$t: \Cal S(\Sigma_{0,4}) \to \bold R_{\geq 2}$ is a trace 
function of a hyperbolic metric if and only if for each ideal triangle $(\alpha_1, 
\alpha_2, \alpha_3)$ so that $(\alpha_i, b_j, b_k)$ bounds 
a $\Sigma_{0,3}$ in $\Sigma_{0,4}$ the following
hold.
$$\prod_{i=1}^3 t(\alpha_i)  = \sum_{i=1}^3 t^2(\alpha_i) + \sum_{j=1}^4
t^2(b_j) + \prod_{j=1}^4 t(b_j) + \frac{1}{2} \sum_{i=1}^3 \sum_{j=1}^4 t(\alpha_i)t(b_j)t(b_k) -4 \quad \text{and} $$
$$ t(\alpha_3)t(\alpha_3') = 
\sum_{i=1}^2 t^2(\alpha_i) + \sum_{j=1}^4
t^2(b_j) + \prod_{j=1}^4 t(b_j) + \frac{1}{2} \sum_{i=1}^2 \sum_{j=1}^4 t(\alpha
_i)t(b_j)t(b_k) -4 $$
where ($\alpha_1$, $\alpha_2$, $\alpha_3'$) and 
 $(\alpha_1, \alpha_2, \alpha_3)$ are two distinct ideal triangles in 
$\Cal S'(\Sigma_{0,4})$.
\rm

Part (a) of theorem A1 was a result of Fricke-Klein [FK] and Keen [Ke2]
and part (b) was proved in [Lu1].

\medskip
\noindent
{\bf Theorem A2.} \it (a) For the surface $\Sigma_{1,1}$, a function $f: \Cal S(
\Sigma_{1,1})$
$\to \bold R_{\geq 0}$  is an intersection  function
if and only if
the following hold. 
$$
f(\alpha_1) + f(\alpha_2) + f(\alpha_3) = \max_{i=1,2,3}(2 f(\alpha_i),
f([\partial \Sigma_{1,1}])) 
$$ where $(\alpha_1, \alpha_2, \alpha_3)$
 is an ideal triangle, and
$$f(\alpha_3) + f(\alpha_3') = \max(2f(\alpha_1), 2f(\alpha_2),
f([\partial \Sigma_{1,1}]))$$ 
where   $(\alpha_1, \alpha_2, \alpha_3)$ and
 $(\alpha_1, \alpha_2, \alpha_3')$  are two distinct ideal triangles.

(b) For the surface $\Sigma_{0,4}$ with $\partial \Sigma_{0,4} = b_1 \cup b_2 \cup
b_3 \cup b_4$, a function $f: \Cal S(\Sigma_{0,4}) \to \bold R_{\geq 0}$ is an
intersection function
if and only if
for each ideal triangle $(\alpha_1, \alpha_2, \alpha_3)$ so that
$(\alpha_i, b_s, b_r)$ bounds a $\Sigma_{0,3}$ in $\Sigma_{0,4}$ the following
hold.
 $$\Sigma_{i=1}^3 f(\alpha_i) = \max_{1 \leq i \leq 3,1 \leq s \leq 4}
(2f(\alpha_i), 2f(b_s), \sum_{s=1}^4 f(b_s), f(\alpha_i)+ f(b_s) + f(b_r))
$$
$$f(\alpha_3) + f(\alpha_3') =\max_{1 \leq i \leq 2,1 \leq s \leq 4}
(2f(\alpha_i), 2f(b_s), \sum_{s=1}^4 f(b_s), f(\alpha_i)+ f(b_s) + f(b_r))
$$
where   $(\alpha_1, \alpha_2, \alpha_3)$ and
$(\alpha_1, \alpha_2, \alpha_3')$  are two distinct ideal triangles.
\rm

Theorem A2 was proved in [Lu2].

Below is the statement of the presentation of the mapping class  group
for all surfaces of negative Euler number (see [Lu3]).

\noindent
{\bf Theorem A3.}  \it For a compact oriented surface $\Sigma$ of
negative Euler number, the
mapping class group $\Gamma(\Sigma)$ has the following presentation:

generators: \rm \{$D_{\alpha} : \alpha \in \Cal S(\Sigma)$\}.

\it relations: \rm 
\text{(I)} $D_{\alpha} D_{\beta} = D_{\beta} D_{\alpha}$ \it
if $\alpha \cap \beta = \emptyset$. \rm

\text{(II)} $D_{\alpha \beta} = D_{\alpha}D_{\beta} D_{\alpha}^{-1}$  \it if
$\alpha \perp \beta$. \rm

\text{(III)} $(D_{\alpha}D_{\beta}D_{\alpha \beta})^4 = D_{\partial(\alpha \cup
\beta)}$ if $\alpha \perp \beta$.
\rm

\text{(IV)} $D_{\alpha}D_{\beta}D_{\alpha \beta} = D_{\partial(\alpha, \beta)}$
\it if $\alpha \perp_0 \beta$. \rm

The characterization of the SL(2, $\bold K)$ characters for surface
group representations is given by the following. The theorem is
proved by Fricke and Klein [FK] and Vogt [Vo], thought stated in different
terminologies.
See for instance [Go] and [Lu3].

\noindent
{\bf Theorem A4.} \it (a) For the  surface $\Sigma_{1,1}$  with $b =\partial \Sigma
_{1,1}$,
a function $t: \Cal S(\Sigma_{1,1}) \to \bold C$  is an $SL(2, \bold C)$ trace
function  if and only if
the following hold.
$$ \prod_{i=1}^3 t(\alpha_i) = \sum_{i=1}^3 t^2(\alpha_i) - t(b) -2 \quad
\text{and}$$
$$ t(\alpha_3) + t(\alpha_3') =  t(\alpha_1) t(\alpha_2)$$
where $(\alpha_1, \alpha_2, \alpha_3)$ and $(\alpha_1, \alpha_2, \alpha_3')$
are distinct ideal triangles in $\Cal S' (\Sigma_{1,1})$.

(b) For the  surface $\Sigma_{0,4}$ with $\partial \Sigma_{0,4} = \cup_{i=1}^4
b_i$,
a function $t: \Cal S(\Sigma_{0,4}) \to \bold C$ is an $SL(2, \bold C)$ trace
function if and only if for each ideal triangle $(\alpha_1,
\alpha_2, \alpha_3)$ so that $(\alpha_i, b_j, b_k)$ bounds
a $\Sigma_{0,3}$ in $\Sigma_{0,4}$ the following
hold.
$$\prod_{i=1}^3 t(\alpha_i)  = - \sum_{i=1}^3 t^2(\alpha_i) - \sum_{j=1}^4
t^2(b_j) - \prod_{j=1}^4 t(b_j) + \frac{1}{2} \sum_{i=1}^3 \sum_{j=1}^4 t(\alpha
_i)t(b_j)t(b_k) +4 \quad \text{and} $$
$$ t(\alpha_3) + t(\alpha_3') = -t(\alpha_1) t(\alpha_2) 
 + \frac{1}{2} \sum_{i=1}^2 \sum_{j=1}^4 t(\alpha _i)t(b_j)t(b_k) $$
where ($\alpha_1$, $\alpha_2$, $\alpha_3'$) and
 $(\alpha_1, \alpha_2, \alpha_3)$ are two distinct ideal triangles in
$\Cal S'(\Sigma_{0,4})$.
\rm

\medskip
\noindent
{\bf Appendix B.} The Modular Structure on the Space of 3-Holed Sphere Decompositions

The other natural example of compact modular structure is the set $HD(\Sigma)$
of all isotopy classes of 3-holed sphere decompositions of a surface $\Sigma$.
The charts are constructed as follows. Suppose $(a_1, ..., a_k)$ is an
element in $HD(\Sigma)$. Take an essential subsurface $\Sigma'$ of
level-1 so that all but one, say $a_i$ of the coordinate, are disjoint from
$\Sigma'$. Now the chart associate to $\Sigma'$ is the set of elements
$\{(a_1, ..., a_{i-1}, b_i, a_{i+1},..., a_k) \in HD(\Sigma) |
b_i \in \Cal S'(\Sigma')\}$ with chart map sending the element to the slope of
$b_i$.  Again if two charts overlap in two elements, they coincide.
A result of Hatcher-Thurston [HT] says that given any two elements
in $HD(\Sigma)$ there is a sequence of charts whose union contains
these two elements so that any two adjacent
charts overlap in at least one  element. On the other hand, each
element in $HD(\Sigma)$ determines a handlebody structure on the
surface $\Sigma$ obtained by attaching 2-cells to the components of
the 3-holed sphere decomposition and then 3-cells. 
Evidently if two elements in $HD(\Sigma)$ lie in a chart 
associated to a 4-holed sphere,
then they determine the same handlebody structure. The main result
in [Lu5] shows that conversely it is also true. Namely,  if two
elements in  $HD(\Sigma)$ determine the same handlebody structure, 
then there  is a  sequence of charts associated to 4-holed spheres
whose union contains these two elements so that any two 
adjacent charts overlap in at least one element.

\medskip
\noindent
{\bf Appendix C.} A Simple Proof That the
Multiplication Produces  a Curve System

For simplicity,
let us assume that the surface is closed (see [Lu2] for general cases).
We will give a simple proof of the following fact used in \S4. Namely,
if $a$ and 
$b$ are two curve systems so that they intersect
minimally in their isotopy classes, then the 1-dimensional submanifold $ab$
obtained by resolving all intersection points in $a \cup b$ from $a$ to
$b$ is again a curve system. 

Suppose otherwise that the 1-submanifold $ab$ contains a null homotopic
component $c$. We may assume that $c$ is the ``inner-most" component, 
i.e., in the interior of the disc  $D$ bounded by $c$, there are no other
components of $ab$.  Let us consider all components of $\Sigma -(a \cup b)$
which are inside $D$, say $A_1,..., A_k$. Each $A_i$ is  an open disc
since $c$ is the inner-most. The boundary of $A_i$ consists of arcs in
$a$ and $b$, and the  corners of $A_i$ corresponds to the intersection points
of $a$ and $b$. Thus
we may call each $A_i$ a polygon bounded by sides in $a$ and $b$
alternatively. Since $a$ intersects $b$ minimally within their
isotopy classes, each $A_i$ has at least four sides. Now by the definition
of the resolution, the disc $D$ is obtained by resolving corners of
$A_i$'s from $a$ to $b$. Considering the resolutions at the vertices  along 
the boundary of $A_i$, one sees that corners open and
close alternatively in a  cyclic order on the boundary.
 Form a graph in $D$ by assigning a vertex
in each $A_i$ and joining an edge between two vertices if their corresponding
polygons $A_i$ and $A_j$ 
have the same corner which is opened by the resolution. Then, on
one hand, the graph is a tree since it is homotopic to the disk.
On the other hand, each vertex of the graph has valence at least two
since the valence of the vertex is half of the number of sides of the
corresponding polygon $A_i$ (by the alternating property). 
This contradicts the fact that a tree must have a vertex of valence one.

\midspace{0.1cm}
\centerline{\epsfbox{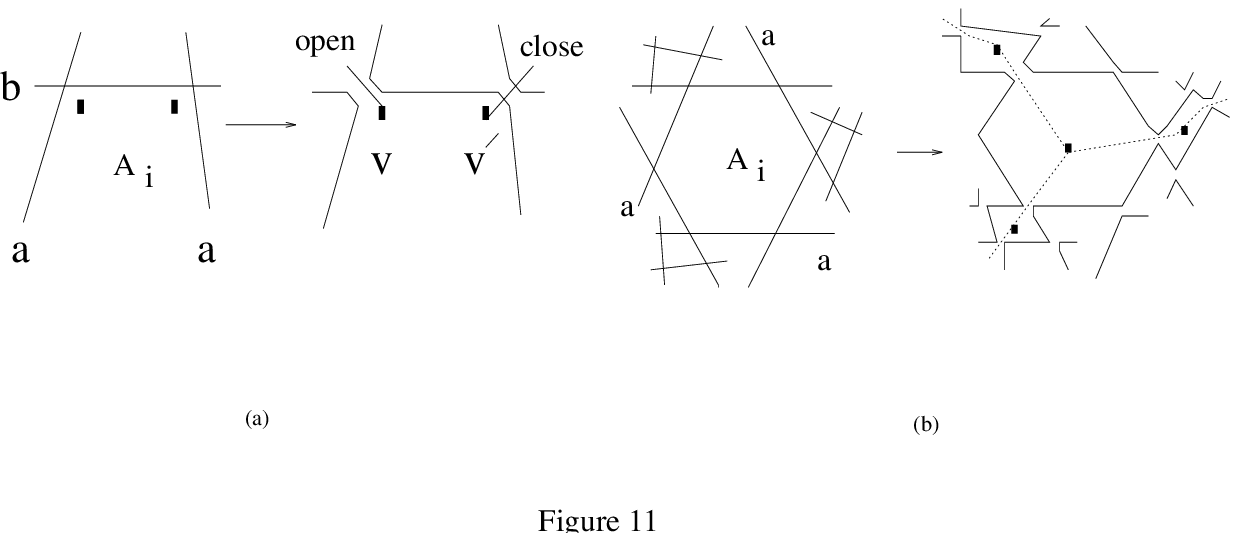}}

\noindent
{\bf Appendix D.} {A Proof of the Reduction Lemma in \S5}

We shall prove the following reduction lemma stated in \S5.  Suppose
$\Sigma$ is a surface of level at least 3. If
$X$ is a subset of
$Y^{\Cal S(\Sigma)}$ so that for each level-2 essential subsurface
$\Sigma'$ the restriction of $X$ to $Y^{\Cal S(\Sigma')}$ has property
RP, then $X$ has property RP. 

To begin the proof, suppose the surface $\Sigma$ is decomposed into a union of two essential
subsurfaces $A_1$ and $A_2$ overlapping in an essential level-0 surface
and we are given two elements $f, g$ in $X$ whose restrictions to
$A_i$ are the same. The goal is to show that $f =g$. To this end,
let us construct a 3-holed sphere decomposition $(a_1, ..., a_k)$ of
the surface $\Sigma$ so that each 3-holed sphere in the decomposition
is either in $A_1$
or in  $A_2$ and $A_1 \cap A_2$ is bounded by $a_i$'s. 
Thus if $s$ is an element in $\Cal S(\Sigma)$ which intersects only
one element of $\{a_1, ..., a_k\}$, then $s$ is in $\Cal S(A_i)$
for $i=1$ or $2$. In particular $f(s) = g(s)$. Now suppose we
make an elementary move on $\{a_1, ..., a_k\}$ to
produce a new 3-holed sphere decomposition $\{b_1, ..., b_k\}$
where all but one of $b_i$ are $a_i$ and the exceptional component,
say $b_j$ is modular related to $a_j$ (i.e., $b_j \perp a_j$
or $b_j \perp_0 a_j$)(these moves were introduced in the appendix
of [HT]).  
We claim that if $f(s) = g(s)$ for all elements $s$ which intersects
at most one of $\{a_1,..., a_k\}$ and $\{b_1,..., b_k\}$
is obtained from $\{a_1,..., a_k\}$ by an elementary move, then
$f(s) = g(s)$ for all elements $s$ which intersects at  most one
of $\{b_1,..., b_k\}$. Indeed, by the 
property RP for level-2 surfaces, we see that $f(s) = g(s)$ for
all elements $s$ inside any level-2 subsurface which is  bounded by
elements in
$\{a_1, ..., a_k\}$.
Now if $s$ is an isotopy class which intersects at most one element
in $\{b_1, ..., b_k\}$, then $s$ intersects at most two elements
in  $\{a_1,..., a_k\}$. Thus the isotopy class $s$ is in a 
level-2 subsurface which is bounded by elements in
$\{a_1,..., a_k\}$. Thus $f(s)= g(s)$.
Now by the result in [HT]  that any two 3-holed sphere decompositions
of the surface are related by a finite sequence of elementary moves,
it follows that $f(s) =g(s)$ for any element $s$ in $\Cal S(\Sigma)$.

In view of the importance of the 3-holed sphere decompositions,
it is attempting to make a 2-dimensional cell-complex $Z$ based on 3-holed sphere decompositions of the surfaces as follows.
The vertices of $Z$ are the isotopy classes of 3-holed sphere decompositions
of the surface and the edges are those pair of vertices related by
an elementary move. Now attaching a 2-cell to each
5-gon associated to each pentagon relation (see figure 9),
a 2-cell to each 4-gon associated to four elementary moves
which are supported in two disjoint level-1 surfaces, and a 2-cell
to each 3-gon associated to three elementary moves supported in
a level-1 surface. 
This cell-complex was implicitly introduced in the appendix of 
[HT]. The simple connectivity of the cell-complex seems to be asserted
in [HT].  See also [Ha] for related topics.

\bigskip
\centerline{\bf Reference}

[Bi] Birman, J.S.,  Braids, links, and mapping class groups.
Ann. of Math. Stud., 82, Princeton Univ. Press, Princeton, NJ, 1975.

[Bo1] Bonahon, F.: Shearing hyperbolic surfaces, bending pleated surfaces and Thurston's symplectic form. Ann. Fac. Sci. Toulouse Math. (6) {\bf 5}
(1996), no. 2, 233--297. 

[Bo2] Bonahon, F.: Bouts des varietes hyperboliques de dimension $3$. 
Ann. of Math. (2) 124 (1986), no. 1, 71--158.

[De] Dehn, M.: Papers on group theory and topology. J. Stillwell (eds.).
Springer-Verlag, Berlin-New York, 1987.

[FK] Fricke, R., Klein, F.: Vorlesungen  \"uber die Theorie der
Automorphen Functionen.  Vol. 2, Teubner, Leipizig, 1897.

[FLP] Fathi, A., Laudenbach, F., Poenaru, V.: Travaux de Thurston sur les
surfaces. Ast\'erisque {\bf 66-67}, Soci\'et\'e Math\'ematique de France, 1979.

[Ge] Gervais, S.: Presentation and central extensions of
mapping class groups. Trans. Amer. Math. Soc. 348 (1996), no. 8, 3097--3132.

[Go] Goldman, W.: Topological components of spaces of representations.
Invent. Math. {\bf 93}(3) (1988), 557-607.

[Gr1] Grothendieck, A.:  Esquisse d'un programme. With an 
English translation by L. Schneps on pp. 243--283. London Math.
Soc. Lecture Note Ser., 242, Geometric Galois actions, 1, 5--48, Cambridge Univ. Press, Cambridge, 1997. 

[Gr2]  Geometric Galois actions. 1. Around Grothendieck's "Esquisse d'un programme". Edited by Leila Schneps and Pierre Lochak. London
Mathematical Society Lecture Note Series, 242. Cambridge University Press, Cambridge, 1997. 

[Gr3] The Grothendieck theory of dessins d'enfants. Papers from the Conference on Dessins d'Enfants held in Luminy, April 19--24, 1993. Edited by
Leila Schneps. London Mathematical Society Lecture Note Series, 200. Cambridge University Press, Cambridge, 1994. 

[Ha] Harer, J.:
Stability of the homology of the mapping class groups of
orientable surfaces. Ann. of Math. (2) 121 (1985), no. 2, 215--249.

[HT] Hatcher, A. Thurston, W.: A presentation for the mapping class group
of a closed orientable surface. Topology {\bf 19} (1980), 221-237.

[Iv] Ivanov, N.: Automorphisms of Teichm\"uller modular groups. In:
O. Ya. Viro (ed.) Topology and geometry. Lecture Notes in Math. Vol.
1346, pp.190-270. Springer-Verlag, Berlin-New York, 1988

[Jo] Johnson, D.:  Homeomorphisms of a surface which act trivially on homology. Proc. Amer. Math. Soc. 75 (1979), no. 1, 119--125.

[Ke1] Keen, L.: private communication.

[Ke2] Keen, L.: Intrinsic moduli on Riemann surfaces. Ann. of Math. {\bf 84} 
(1966), 405-420.

[Li] Lickorish, R.: A representation of oriented combinatorial 3-manifolds. 
Ann. of
Math. {\bf 72} (1962), 531-540.

[Lu1] Luo, F.:  Geodesic length functions and Teichm\"uller spaces. \it
J. Differential Geom. \rm  {\bf 48}, (1998), 275-317.

[Lu2] Luo, F.: Simple loops on surfaces and their intersection numbers, 
preprint, 1997. The Research announcement appeared in  \it Math. Res. Letters, \rm {\bf 5}, (1998), 47-56.

[Lu3] Luo, F.: SL(2, $\bold C$) characters for surface groups, preprint, 1999.

[Lu4] Luo, F.: A presentation of the mapping class group.
\it Math. Res. Letters, \rm {\bf 4}, (1997), 735-739.

[Lu5] Luo, F.:  Automorphisms of the curve complex, \it Topology, \rm
to appear.

[Lu6] Luo, F.:  On Heegaard diagrams. \it Math. Res. Letters,
\rm  {\bf 4}, (1997), 365-373.

[Lu7] Luo, F.: Singular affine flat sturctures on surfaces, in 
preparation.

[Ma] Maskit, B.: A picture of moduli space. Invent. Math. 126 (1996), no. 2,
341-390. 

[MS] Moore, C. and Seiberg, N.: Polynomial equations for rational
conformal field theories, Phys. Lett. B. {\bf 212}, (1988), 451-460.

[Mo] Mosher, L.: Tiling the projective foliation space of a punctured surface.
Trans. Amer. Math. Soc.  {\bf 306}, (1988), 1-70.

[Ok1] Okumura, Y.: On the global real analytic coordinates for Teichm\"uller
spaces. J. Math. Soc. Jap. {\bf 42} (1990), 91-101

[Ok2] Okumura, Y.: Global real analytic length parameters for
Teichm\"uller spaces. Hiroshima Math. J.  {\bf 26} (1) (1996), 165--179

[Pe] Penner, R.: The decorated Teichm\"uller space of punctured surf
aces. Comm. Math. Phys. {\bf 113}, (1987), no. 2, 299--339.

[Re] Rees, M.: An alternative approach to the ergodic theory of
measured foliations on surfaces. Ergodic Theory Dynamical Systems 1 (1981),
no. 4, 461-488.

[Sc] Schmutz, P.: Die Parametrisierung des Teichm\"ullerraumes durch
geod\"atishe L\"angenfun \newline ktionen. Comm. Math. Hel. {\bf 68} (1993), 278
-288

[Se] Series, C.: The modular surface and continued fractions,
J. London Math. Soc. {\bf 31} (1985) no. 2, 69-80.

[Th1] Thurston, W.: Geometry and topology of 3-manifolds, Princeton University
lecture notes, 1979.

[Th2] Thurston, W.: On the geometry and dynamics of diffeomorphisms of
surfaces. Bul. Amer. Math. Soc.  {\bf 19} (1988) no 2, 417-438.       

[Vo] Vogt, H.: Sur les invariants fondamentaux des \'equations 
diff\'erentielles lin\'eaires du second ordre. Ann. de l'Ecole Normale Superieur (3), 6 Suppl.
(1889), 3-72.

[Wo1] Wolpert, S.: Fenchel-Nielsen deformation. Ann. of Math. (2) {\bf 115}, 
(1982), no. 3, 501--528. 

[Wo2] Wolpert, S.: private communication.

Department of Mathematics, Rutgers University, New Brunswick, NJ 08903 

\end